\documentclass[10pt]{amsart} 

\usepackage[utf8]{inputenc}
\usepackage[T1]{fontenc}	
\usepackage[french,english]{babel}	

\usepackage{lmodern}			

\usepackage{graphicx}	

\usepackage[a4paper,plainpages=false,colorlinks,linkcolor=bleuFonce,citecolor=rougeFonce,urlcolor=vertFonce,breaklinks]{hyperref}

\usepackage{placeins}
\usepackage{float} 

\usepackage{color}
\definecolor{vertFonce}{rgb}{0,0.5,0}
\definecolor{numLignes}{rgb}{0.17,0.57,0.7}	
\definecolor{gris}{rgb}{0.5,0.5,0.5}
\definecolor{grisFonce}{rgb}{0.2,0.2,0.2}
\definecolor{orange}{rgb}{1,0.65,0.31}		
\definecolor{orangeFonce}{rgb}{1,0.4,0}
\definecolor{bleuFonce}{rgb}{0,0,0.4}
\definecolor{rougeFonce}{rgb}{0.3,0,0}
\definecolor{rougeWord}{rgb}{0.5,0,0}
\definecolor{vertClair}{rgb}{0.8,1,0.8}
\definecolor{rougeClair}{rgb}{1,0.5,0.5}

\usepackage{pict2e}
\usepackage{multido}

\usepackage{amsfonts,amssymb,amsthm,amsmath} 
\usepackage{dsfont}				
\usepackage{mathrsfs}
\usepackage{yfonts}
 
\newtheorem{theorem}{Theorem}
\newtheorem{cor}{Corollary}[section]
\newtheorem{prop}{Proposition}[section]

\newtheorem{remark}{Remark}[section]

\newenvironment{demo}[1][]{%
	\begin{proof}[\textbf{Proof#1}]
	}{%
	\end{proof}
}

\newenvironment{thm}[1][]{%
	\color{rougeWord}\begin{theorem}[#1]
	}{%
	\end{theorem}
}

\newcommand{\step}[1]	{\paragraph{\itshape\bfseries Step #1.}}

\usepackage{stmaryrd}
\newcommand		{\subsetArrow}	{\mathrel{\ooalign{$\subset$\cr%
\hidewidth\raise-.087ex\hbox{$_\shortrightarrow\mkern-1.5mu$}\cr}}}
\newcommand		{\subsetarrow}	{\mathrel{\ooalign{$\subset$\cr%
\hidewidth\raise-1.45ex\hbox{$\vec{}\mkern6mu$}\cr}}}

\newcommand		{\N}			{\mathbb N}			
\newcommand		{\R}			{\mathbb R}			
\newcommand		{\CC}			{\mathbb C}			
\renewcommand	{\L}			{\mathcal L}		
\newcommand		{\F}			{\mathcal F}		
\newcommand		{\B}			{\mathscr B}		
\renewcommand	{\P}			{\mathcal P}		
\newcommand		{\PP}			{\mathscr P}		

\newcommand		{\bb}			{\mathfrak b}

\newcommand		{\ssi}			{\Leftrightarrow}

\newcommand		{\il}			{[\![}				
\newcommand		{\ir}			{]\!]}				
\newcommand		{\Int}[1]		{\il #1 \ir}

\newcommand		{\lt}			{\left}
\newcommand		{\rt}			{\right}

\newcommand		{\ii}			{\mathrm{i}}	
\newcommand		{\init}			{\mathrm{in}}
\newcommand		{\loc}			{\mathrm{loc}}

\renewcommand	{\d}			{\,\mathrm{d}}	
\newcommand		{\dt}			{\frac{\mathrm{d}}{\mathrm{d}t}}	
\newcommand		{\ddt}[1]		{\frac{\mathrm{d}#1}{\mathrm{d}t}}	
\DeclareMathOperator{\divg}		{div}

\DeclareMathOperator*{\supess}	{sup\,ess}
\DeclareMathOperator{\diag}		{diag}

\newcommand		{\intd}			{\int_{\R^d}}
\newcommand		{\intdd}		{\int_{\R^{2d}}}
\newcommand		{\iintd}		{\iint_{\R^{2d}}}
\newcommand		{\sumj}			{\sum_{j\in J}}
\newcommand		{\eps}			{\varepsilon}

\newcommand		{\Eps}			{\mathcal{E}}

\usepackage{upgreek}
\renewcommand	{\r}		{\uprho}			
\newcommand		{\op}		{\boldsymbol{\rho}}	
\newcommand		{\opmu}		{\boldsymbol{\mu}}	
\newcommand		{\gam}		{\boldsymbol{\gamma}}
\DeclareMathOperator{\OP}	{op}
\DeclareMathOperator{\Tr}	{Tr}				
\newcommand		{\conj}[1]	{\overline{#1}}		
\newcommand		{\bra}[1]	{\langle #1 |}
\newcommand		{\ket}[1]	{| #1 \rangle}
\newcommand		{\Wh}		{W_{2,\hbar}}		
\newcommand		{\wt}		{\widetilde{w}}		
\newcommand		{\w}		{w}					
\newcommand		{\wh}		{\w_\hbar}					

\newcommand		{\opp}		{\boldsymbol{p}}

\newcommand		{\gh}		{\gam_\hbar}
\newcommand		{\ch}		{\mathbf{c}_\hbar}
\newcommand		{\Eh}		{\mathcal{E}_\hbar}
\newcommand		{\fh}		{f_\hbar}
\newcommand		{\fht}		{\tilde{\fh}}
\newcommand		{\oph}		{\op_\hbar}
\newcommand		{\Hh}		{H_\hbar}

\title[Propagation of Moments and Limit from Hartree to Vlasov Equation]{Propagation of Moments and Semiclassical Limit from Hartree to Vlasov Equation}
\author{Laurent Lafleche}
\date{\today}

\def\signll{\bigskip\begin{center}{
	\sc Laurent
	Lafleche\par\vspace{3mm} 
	Université Paris-Dauphine, PSL\par
	CEREMADE, UMR CNRS 7534\par
	Place du Maréchal de Lattre de Tassigny \par
	75775 Paris Cedex 16 FRANCE\par\vspace{3mm} e-mail:}
	\tt{lafleche@ceremade.dauphine.fr}
\end{center}}

\begin{document}

\begin{abstract}In this paper, we prove a quantitative version of the semiclassical limit from the Hartree to the Vlasov equation with singular interaction, including the Coulomb potential. To reach this objective, we also prove the propagation of velocity moments and weighted Schatten norms which implies the boundedness of the space density of particles uniformly in the Planck constant.
\end{abstract}

\maketitle

\bigskip

\renewcommand{\contentsname}{\centerline{Table of Contents}}
\setcounter{tocdepth}{2}	
\tableofcontents


\bigskip
\section{Introduction}\label{sec:intro}

\subsection{Presentation of the problem}

	In this paper, we consider the nonrelativistic quantum and classical equations which describe the evolution of a density of infinitely many particles in the kinetic mean field regime, called respectively the Vlasov and the Hartree equation. The interaction between particles is described by a mean field potential $V = V(x)$ depending only on the space variable $x\in\R^d$ with $d\geq 2$ and which is defined by
	\begin{equation*}
		V := K * \rho = \intd K(x-y)\rho(y)\d y,
	\end{equation*}
	where $\rho$ is the spatial density and $K$ is an even kernel describing the interaction between two particles. The force field can then be written 
	\begin{equation*}
		E := -\nabla V.
	\end{equation*}
	Typically, we have in mind the pair interaction potential $K(x) = \frac{\pm 1}{|x|^a}$ with $a\in [-2,d-1)$. The most physically relevant case is the case of the Coulomb interaction $a=d-2$ for $d\geq 3$ or $K(x) = \pm\ln(|x|)$ in the two dimensional case. It can describe the interaction of charged particles as well as a system of point masses in gravitational interaction, the force being repulsive when $K$ is positive and attractive in the converse case.
	
	In the classical case, the kinetic density of particles $f = f(t,x,\xi)$ is a nonnegative function of time $t\in\R_+$, space and momentum $\xi\in\R^d$ and the space density is given by
	\begin{align*}
		\rho(x) &:= \intd f(t,x,\xi)\d \xi.
	\end{align*}
	The evolution of the kinetic density is then given by the well-known Vlasov equation
	\begin{align}\tag{Vlasov}\label{eq:Vlasov}
		\partial_t f + \xi\cdot\nabla_xf + E\cdot\nabla_\xi f & = 0.
	\end{align}
	Remark also that by defining the Hamiltonian
	\begin{align*}
		H &:= \frac{|\xi|^2}{2} + V,
	\end{align*}
	we can write the \eqref{eq:Vlasov} equation as
	\begin{align*}
		\partial_t f &= \lt\{H,f\rt\},
	\end{align*}
	where $\{\cdot,\cdot\}$ is the Poisson bracket defined by
	\begin{equation*}
		\{u,v\} = \nabla_x u\cdot\nabla_\xi v - \nabla_\xi u\cdot\nabla_x v.
	\end{equation*}
	
	On the other hand, in the formalism of quantum mechanics, a particle is described by a wave function $\psi \in L^2 = L^2(\R^d,\CC)$ verifying $\|\psi\|_{L^2}=1$. Under the action of the potential $V$, its evolution is governed by the following Schrödinger equation
	\begin{align}\label{eq:NLS}
		i\hbar\partial_t\psi & = -\frac{\hbar^2}{2}\Delta \psi + V\psi,
	\end{align}
	where $\hbar = \frac{h}{2\pi}$ is the reduced Planck constant. In the more general case of systems with mixed states, the density of particles is described by a trace class and self-adjoint density operator, $\op$, which by the Spectral theorem can be seen as a superposition of pure orthonormal states $(\psi_j)_{j\in J} \in (L^2)^J$ for a given $J\subset\N$ by writing
	\begin{align}\label{def:op}
		\op \varphi &:= \intd \r(x,y)\varphi(y)\d y = \sumj \lambda_j \ket{\psi_j}\bra{\psi_j}.
	\end{align}
	This is a Hilbert-Schmidt operator of kernel
	\begin{equation*}
		\r(x,y) = \sumj \lambda_j \conj{\psi_j(y)}\psi_j(x).
	\end{equation*}
	Given a density operator, the spatial density is defined as the diagonal of the kernel 
	\begin{align*}	
		\rho(x) &:= \r(x,x) = \sumj \lambda_j |\psi_j|^2,
	\end{align*}
	and the Hamiltonian is the following operator
	\begin{align*}
		H &= -\frac{\hbar^2}{2}\Delta + V
	\end{align*}
	where $V = K*\rho$ is identified with the operator of multiplication by $V(x)$. We can then rewrite \eqref{eq:NLS} for each $\psi_j$ as $\partial_t\psi_j = \frac{1}{i\hbar}H\psi_j$ and we deduce that the density operator verifies the so called Hartree equation
	\begin{align}\label{eq:Hartree}\tag{Hartree}
		\partial_t\op &= \frac{1}{i\hbar}\lt[H,\op\rt],
	\end{align}
	where $[\cdot,\cdot]$ is the Lie bracket defined by
	\begin{equation*}
		[A,B] = AB-BA.
	\end{equation*}
	
	The main goal of the present paper is to obtain a \textit{quantitative estimate} of the semiclassical limit from the \eqref{eq:Hartree} to the \eqref{eq:Vlasov} equation, which means the limit when $h = 2\pi\hbar \to 0$. This limit was first investigated in a non-quantitative way  using compactness methods by Lions and Paul \cite{lions_sur_1993}, Markowich and Mauser \cite{markowich_classical_1993} and then by Gerard et al \cite{gerard_homogenization_1997}, Gasser et al \cite{gasser_semiclassical_1998}, Ambrosio et al \cite{ambrosio_passage_2010, ambrosio_semiclassical_2011}, Graffi et al \cite{graffi_mean-field_2003}. On the other hand, Athanassoulis et al \cite{athanassoulis_strong_2011} prove quantitative estimates in $L^2$ norm in the case of sufficiently smooth potentials and Amour et al \cite{amour_classical_2013, amour_semiclassical_2013} show that the rate can be improved in the case of very smooth potentials. More recently, some improvements on the requirement of regularity of the potential $K$ have been done in Benedikter et al \cite{benedikter_hartree_2016} by considering trace and Hilbert-Schmidt norms and a mixed semiclassical and mean-field limit, and by Golse et al \cite{golse_mean_2016} and Golse and Paul \cite{golse_schrodinger_2017} using quantum pseudo-distances created on the model of the Wasserstein-Monge-Kantorovitch distances. This strategy allows them to prove estimates that do not require any assumption of regularity on the initial data. However, all these works still require at least Lipschitz regularity of the potential, which does not include singular interactions like the Coulomb potential.
	
	Recent attempts on generalizing these results to more singular potentials in the case of fermionic systems can be found in the works by Porta et al \cite{porta_mean_2017} and Saffirio \cite{saffirio_mean-field_2018}, where a joint mean-field and semiclassical limit is obtained. However, it requires regularity assumptions on the solution of the Hartree equation whose propagation is still an open problem. The closely related problem of the mean field limit from the $N$-body Schrödinger to the Hartree equation has been also investigated a lot. Weak convergence results have been first obtained in \cite{bardos_weak_2000, erdos_derivation_2001, bardos_derivation_2002} for the Coulomb potential. See also \cite{zhang_limit_2002} for the one dimensional case. Quantitative results have been established in \cite{rodnianski_quantum_2009, pickl_simple_2011, golse_mean_2016, mitrouskas_bogoliubov_2016, golse_schrodinger_2017, golse_empirical_2017, golse_derivation_2018} for Bosons and in \cite{frohlich_mean-field_2009, frohlich_microscopic_2011, benedikter_mean-field_2014, benedikter_mean-field_2016, bach_kinetic_2016, petrat_new_2016, porta_mean_2017, petrat_hartree_2017} for Fermions.	Remark that some of these works use a joint mean-field and semiclassical limit, however they always require at least a Lipschitz potential or an assumption of regularity on the solution of the Hartree equation.
	
	An other possible way to derive the Vlasov equation is the classical mean-field limit. It is also a closely related problem. Results for non-smooth potentials can be found for example in 	\cite{hauray_n-particles_2007, hauray_particle_2015, lazarovici_vlasov-poisson_2016, jabin_mean_2016, lazarovici_mean_2017}. The major obstacle here is the absence of regularity in the $N$-body problem, which is the reason why all results with unbounded pair interaction potentials need a cut-off on the force field.
	
	Other results about the mean-field limit are the convergence of the minimizers of the $N$-particles energy towards the mean-field energy. We refer for example to \cite{fournais_semi-classical_2018} and references therein.
	
	In order to get semiclassical estimates for more general pair potentials $K$, our strategy consists in requiring more regularity on the initial data and proving that it implies regularity at the level of the mean-field potential $V$. The propagation of moments is inspired from 	\cite{lions_propagation_1991} and \cite{lions_sur_1993}. The semiclassical limit is mostly an adaptation of \cite{golse_schrodinger_2017} and of the proof of uniqueness for the Vlasov equation given in \cite{loeper_uniqueness_2006}. Some interesting improvements for the uniqueness can be found in \cite{miot_uniqueness_2016, holding_uniqueness_2018}.

	Finally, notice that the global well-posedness in Sobolev and Schatten spaces and conservation of the energy have been treated in \cite{ginibre_class_1980, ginibre_global_1985, hayashi_smoothing_1989, brezzi_three-dimensional_1991, illner_global_1994, castella_l2_1997, lewin_hartree_2015}. In particular, our hypotheses on finite quantum moments of order $n$ require the equation to be well-posed in the corresponding  $H^n$ Sobolev space. However, as we will see, even if quantum moments can be interpreted as $H^n$ norms, the above mentioned papers do not prove the propagation of these norms uniformly with respect to $\hbar$, which is one the main results of this paper.

\subsection{Notations and tools}

	We describe in this section the main notations that we will use. Since we are in the semiclassical regime, most of our results have to be true in the limit and are inspired from the classical results. Therefore, our notations try to be close for the classical objects and their quantum counterpart. When comparing quantum and classical objects, we will sometimes add $\hbar$ in the notation to denote the quantum objects.

	\subsubsection{Functional spaces}
	Since most of the functional spaces we use will be defined for functions defined on $\R^d$, we will often write $X = X(\R^d,\CC)$, as for example in the case of the Lebesgue spaces $L^p := L^p(\R^d,\CC)$. When working on the phase space $\{(x,\xi)\in \R^{2d}\}$ we will write $L^p_{x,\xi} := L^p(\R^{2d},\CC)$. Some other standard functional spaces we will use are the weak and weighted Lebesgue spaces, defined reciprocally by
	\begin{align*}
		L^{p,\infty} &:= \{f \text{ measurable}, \forall\lambda>0,|\{|f|>\lambda\}|\leq C/\lambda^p\}
		\\
		L^{\infty,\infty} &:= L^\infty
		\\
		L^p(m) &:= \{f \text{ measurable}, fm\in L^p\}.
	\end{align*}
	Moreover, we will denote by $\P(X)$ the space of probability measures on some space $X$. We will need the equivalent of some of these spaces in the quantum picture. The quantum equivalent of the integral on the phase space is the trace which for an operator $\op$ in the form \eqref{def:op} can be written
	\begin{equation*}
		\Tr(\op) = \intd \r(x,x)\d x = \sum_{j\in J} \lambda_j.
	\end{equation*}
	The trace is defined more generally for trace class operators. We refer to \cite{simon_trace_2005} for the general definition and additional properties. In order to define the equivalent of Lebesgue norms, let us first recall the definition of the Schatten norm of a trace class operator $A$ for $p\in[1,+\infty)$
	\begin{align*}
		\|A\|_{p} &:= \Tr(|A|^p)^{1/p}
		\\
		\|A\|_{\infty} &:= \|A\|_{\B},
	\end{align*}
	where $\B = \B(L^2)$ is the space of bounded operator on $L^2$ and $|A| = \sqrt{A^*A}$. We will more precisely use a rescaled version of these norms defined for $r\in[1,+\infty]$ by
	\begin{align}\label{def:q_schatten}
		\|\op\|_{\L^r} &:= h^{-d/r'}\|\op\|_r,
	\end{align}
	where $r' = \frac{r}{r-1}$ denotes the Hölder conjugate of $r$. They play the role of the $L^r_{x,\xi}$ norm for the quantum density operators. The space of quantum probability measures corresponds to the space of normalized hermitian operators defined by 
	\begin{equation*}
		\PP := \{\op\in\B(L^2), \op = \op^* \geq 0, \Tr(\op)= 1\}.
	\end{equation*}
	Remark that since $\op$ will usually be a nice compact operator, for general unbounded operators $A\in\L(L^2)$, we can define $\Tr(A\op) := \Tr(\op^{1/2}A\op^{1/2})$ even if $A\op$ is not a bounded operator.
	
	\subsubsection{Momentum}
	
	We recall that the quantum equivalent of the classical momentum $\xi$ is the following unbounded operator from $L^2$ to $(L^2)^d$
	\begin{align*}
		\opp &:= -i\hbar\nabla.
	\end{align*}
	Its formal adjoint for the scalar product defined by $\langle u,v\rangle_{(L^2)^d} = \intd u\cdot v$ is then defined by
	\begin{align*}
		\opp^* &= -i\hbar\divg = \opp\cdot,
	\end{align*}
	which leads to the following notations
	\begin{align*}
		-\hbar^2\Delta &= |\opp|^2
		\\
		H &= \frac{|\opp|^2}{2} + V.
	\end{align*}
	
	\subsubsection{Wigner Transform}
	There exists several ways to try to associate a density over the phase space to a density operator, one of them being the Wigner transform and its nonnegative but smoothed version called Husimi transform defined reciprocally for $h=1$ by
	\begin{align*}
		\w(\op)(x,\xi) &:= \intd e^{-2i\pi y\cdot\xi}\r\lt(x+\frac{y}{2},x-\frac{y}{2}\rt)\d y = \F(\tilde{\r}_x)(\xi)
		\\
		\wt(\op) &:= \w(\op) * G,
	\end{align*}
	where $\tilde{\r}_x(y) = \r(x+y/2,x-y/2)$ and $G(z) = \frac{1}{\pi^d}e^{-|z|^2}$ with $z := (x,\xi)$ and we used the following convention for the Fourier transform
	\begin{align*}
		\F(u)(\xi) := \intd e^{-2i\pi x\cdot\xi}u(x)\d x.
	\end{align*}
	We refer for example to \cite{lions_sur_1993} and \cite{golse_mean_2016} for more details and mathematical results. Given $\op$ solution of the $\eqref{eq:Hartree}$ equation we will write its Wigner and Husimi transforms respectively
	\begin{align*}
		\fh(x,\xi) &= \wh(\op)(x,\xi) := \frac{1}{h^d}\w(\op)\lt(x,\frac{\xi}{h}\rt)
		\\
		\fht(x,\xi) &:= \fh * G_\hbar,
	\end{align*}
	where $G_\hbar(z) = \frac{1}{(\pi\hbar)^d}e^{-|z|^2/\hbar} = g_\hbar(x)g_\hbar(y)$ with $g_\hbar(x) = \frac{1}{(\pi\hbar)^{d/2}}e^{-|x|^2/\hbar}$. We also define the quantum velocity moments by
	\begin{equation*}
		M_n := \Tr(|\opp|^n\op) = \intdd \fh |\xi|^n\d x\d\xi.
	\end{equation*}
	Remark that the scaling of the quantum Lebesgue norm $\L^p$ can be understood by looking at the Wigner transform and noticing that when $r=1$ or $r=2$
	\begin{align*}
		\iintd \fh\d x\d \xi &= \|\oph\|_{\L^1}
		\\
		\|\fh\|_{L^2_{x,\xi}} &= \|\oph\|_{\L^2}.
	\end{align*}
	Moreover, when $r>2$ and $\oph$ is a superposition of coherent states, then
	\begin{align*}
		\|\fh\|_{L^r_{x,\xi}} &\leq  \|\oph\|_{\L^r}
		\\
		\|\fh\|_{L^r_{x,\xi}} &\underset{h\to 0}{\rightarrow}  \|\oph\|_{\L^r}.
	\end{align*}
	See Section~\ref{sec:coh_stat} for the proof and other results for coherent states.

	\subsubsection{Semiclassical Wasserstein pseudo-distances.} A last useful tool in the study of uniqueness and stability estimates for the Vlasov equation is the Wasserstein-(Monge-Kantorovich) distance $W_p$ which can be defined for any $p\in[1,\infty]$. We refer for example to the books by Villani \cite{villani_topics_2003} and Santambrogio \cite{santambrogio_optimal_2015}. As introduced in \cite{golse_schrodinger_2017}, we will use a quantum equivalent of the $W_2$ distance. We first introduce the notion of coupling between a density operator and a classical kinetic density. Let $\gam\in L^1(\R^{2d},\PP)$. We say that $\gam$ is a semiclassical coupling of $f\in L^1\cap\P(\R^{2d})$ and $\op\in\PP$ and we write $\gam\in\mathcal{C}(f,\op)$ when 
	\begin{align*}
		\Tr(\gam(z)) &= f(z)
		\\
		\intdd \gam(z)\d z &= \op.
	\end{align*}
	Then we define the semiclassical Wasserstein-(Monge-Kantorovich) pseudo-distance in the following way
	\begin{equation}\label{def:Wh}
		\Wh(f,\op) := \lt(\inf_{\gam\in\mathcal{C}(f,\op)}\intdd \Tr\lt(\mathbf{c}_\hbar(z)\gam(z)\rt)\d z\rt)^\frac{1}{2},
	\end{equation}
	where $\mathbf{c}_\hbar(z)\varphi(y) = \lt(|x-y|^2 + |\xi-\opp|^2\rt)\varphi(y)$, $z = (x,\mathbf{\xi})$ and $\opp = -i\hbar\nabla_y$. This is not a distance but it is comparable to the classical Wasserstein distance $W_2$ between the Wigner transform of the quantum density operator and the normal kinetic density, in the sense of the following Theorem
	\begin{theorem}[Golse \& Paul \cite{golse_schrodinger_2017}]\label{th:comparaison_W2}
		Let $\op\in\PP$ and $f\in\P(\R^{2d})$ be such that 
		\begin{equation*}
			\intdd f(x,\xi)(|x|^2+|\xi|^2)\d x\d\xi < \infty.
		\end{equation*}
		Then one has $\Wh(f,\op)^2 \geq d\hbar$ and for the Husimi transform $\fht$ of $\op$, it holds 
		\begin{equation}\label{eq:comparaison_W2}
			W_2(f,\fht)^2 \leq \Wh(f,\op)^2 + d\hbar.
		\end{equation}
	\end{theorem}
	
	See also \cite{golse_wave_2018} for more results about this pseudo-distance and Section~\ref{sec:coh_stat} for the particular case of coherent states.

\subsection{Main results}
	
	We will use this pseudo-distance to get explicit speed of convergence in $\hbar$ of the solution $\oph$ of \eqref{eq:Hartree} equation to the solution $f$ of \eqref{eq:Vlasov} equation. For the classical density $f$, we consider conditions which ensure existence and uniqueness of the solution and the boundedness of $\rho$ as claims the following theorem
	\begin{theorem}[Lions \& Perthame \cite{lions_propagation_1991}, Loeper \cite{loeper_uniqueness_2006}]\label{th:VP_classique}
		Assume $f^\init\in L^\infty_{x,\xi}(\R^6)$ verify
		\begin{equation}\label{hyp:VP_moment}
			\int_{\R^6} f^\init|\xi|^{n_0}\d x\d\xi < C \text{ for a given } n_0 > 6,
		\end{equation}
		and for all $R>0$,
		\begin{equation}\label{hyp:VP_intg}
			\supess_{(y,w)\in\R^6}\{f^\init(y+t\xi,w), |x-y|\leq Rt^2, |\xi-w|\leq Rt\} \in L^\infty_\loc(\R_+,L^\infty_xL^1_\xi).
		\end{equation}
		Then there exists a unique solution to the \eqref{eq:Vlasov} equation with initial condition $f_{t=0} = f^\init$. Moreover, in this case, the spatial density verifies
		\begin{equation}
			\rho\in L^\infty_\loc(\R_+,L^\infty).
		\end{equation}
	\end{theorem}
	
	This is actually proved for the Vlasov-Poisson equation only (i.e. $K = \frac{1}{|x|}$) but the proof would works for less singular potentials verifying the assumptions of the following Theorem~\ref{th:CV}. Actually, the proof we make for the quantum case can be easily adapted to the classical case, which implies for example that this result holds in dimension $3$ for
	\begin{equation}
		K = \frac{1}{|x|^a} \text{ for any } a\in(-1,4/5),
	\end{equation}
	and for all $t\in[0,T_\mathrm{max}]$ when $a\in[4/5,8/7)$. The strategy to prove the above theorem is to obtain a Gronwall's inequality for moments. Our first Theorem uses the same strategy in the semiclassical picture to prove the propagation of quantum velocity moments.
	
	\begin{thm}\label{th:propag_moments}
		Let $r\in [1,\infty]$, $\bb_n := \frac{nr'+d}{n+1}$ and $\nabla K\in L^{\bb,\infty}$ for a given $\bb\in(\max(\bb_4,\bb_n),+\infty]$ and assume $\oph$ verify the \eqref{eq:Hartree} equation for $t\in[0,T]$ with initial condition $\oph^\init\in\PP\cap\L^r$ such that $M_n^\init$ is bounded independently of $\hbar$ for a given $n\in 2\N$. Then there exists $T>0$ and $\Phi\in C^0[0,T)$ such that for any $t\in[0,T)$
		\begin{equation}
			M_n \leq \Phi(t).
		\end{equation}
		Moreover, $T=+\infty$ when $\bb \geq \bb_2 = \frac{2r' + d}{3}$. In particular, if $K=\frac{1}{|x|^a}\in L^{d/a,\infty}$ and $r=\infty$, then we require
		\begin{itemize}
		\item $a \in \lt(-1,\frac{2}{3}\rt)$ if $d=2$,
		\item $a \in \lt(-1,\frac{8}{7}\rt)$ if $d=3$,
		\end{itemize}
		and $T=+\infty$ when
		\begin{itemize}
		\item $a \in \lt(-1,\frac{1}{2}\rt)$ if $d=2$,
		\item $a \in \lt(-1,\frac{4}{5}\rt)$ if $d=3$.
		\end{itemize}
	\end{thm} 
	
	From the quantum kinetic interpolation inequalities \eqref{eq:Lieb_Thirring_+}, we obtain the following corollary.
	\begin{cor}
		Under the assumptions of Theorem~\ref{th:propag_moments},
		\begin{equation*}
			\|\rho\|_{L^p},
		\end{equation*}
		is bounded on $[0,T)$ independently of $\hbar$ for any $p\in[1,p_n]$, where $p_n' = r' + \frac{d}{n}$.
	\end{cor}
	
	\begin{remark}
		As it can be seen in the proof, when $\bb \geq \max(\bb_n,\bb_N)$ for $4\leq n<N\in 2\N$ and there is propagation of moments of order $n$, then the propagation of all higher moments $M_N$ holds with $T=+\infty$. In particular, in dimension $d=3$, if $r=\infty$, as long as the moments $M_4$ are finite, then all the higher moments are propagated for $\bb \geq \frac{7}{5}$, or equivalently if $K=|x|^{-a}$, for $a\leq \frac{8}{7}$, which includes the Coulomb case.
	\end{remark}

	\begin{remark}
		As explained in Section~\ref{sec:energy}, the constraint $a>-1$ could be easily removed by assuming bounded space moments $N_k = \Tr(|x|^k\op)$ for a given $k\leq n$, allowing for polynomial growth of $K$ for large $|x|$.
	\end{remark}
	
	The next theorem is about the following semiclassical convergence result which uses only hypothesis on initial velocity moments and quantum Schatten norms.
	
	\begin{thm}\label{th:CV}
		Let $r\geq2$, $\bb_n := \frac{nr'+d}{n+1}$ and assume $K$ verifies
		\begin{align}\label{hyp:grad_1}
			\nabla K &\in L^\infty + L^{\bb,\infty} &\text{ for some } &\bb\in(\bb_4,+\infty)
			\\\label{hyp:grad_2}
			\nabla^2 K &\in L^2 + L^q &\text{ for some } &q\in(r',2),
		\end{align}
		and let $f$ be a solution of the \eqref{eq:Vlasov} equation and $\oph$ be a solution of \eqref{eq:Hartree} equation with respective initial conditions 
		\begin{align*}
			f^\init &\in \P\cap L^\infty_{x,\xi} \text{ verifying \eqref{hyp:VP_moment} and \eqref{hyp:VP_intg}}
			\\
			\oph^\init &\in \PP\cap \L^r.
		\end{align*}
		Assume also that the initial quantum velocity moment
		\begin{equation}
			M_{n_1}^\init < C \text{ for a given } n_1 \geq \frac{d}{q-r'}.
		\end{equation}
		Then there exists $T>0$ such that for any $t\in(0,T)$,
		\begin{equation*}
			\Wh(f(t),\oph(t)) \leq C_T \lt(\Wh(f^\init,\oph^\init) + \sqrt{\hbar}\rt).
		\end{equation*}
		Moreover, when $\bb \geq \frac{d+2r'}{3}$, then there exists $\Phi\in C^0(\R_+)$ such that for any $t>0$
		\begin{equation}\label{eq:CV_t}
			\Wh(f(t),\oph(t)) \leq \Wh(f^\init,\oph^\init) e^{C(t)} + C_0(t)\sqrt{\hbar},
		\end{equation}
		where
		\begin{align*}
			C_1(t) &= \|\nabla^2 K\|_{L^{s,\infty}}^2 \Phi(t)^2
			\\
			C(t) &= 1 + C_1(t) + \|\nabla^2 K\|_{L^{q}}\Phi(t)
			\\
			C_0(t) &= C_1(t)C(t)^{-1}(e^{2C(t)}-1).
		\end{align*}
	\end{thm}
	
	 The next theorem proves the semi-classical convergence in a case of more singular interactions kernels such that $\nabla K$ is in the Besov space $B^1_{1,\infty}$, which includes the Coulomb potential. The definition and basic properties of Besov spaces are recalled in Appendix~\ref{appendix:Besov}.
	
	\begin{thm}\label{th:CV_VP}
		Assume $K$ verifies
		\begin{align*}
			\nabla K &\in L^\bb + L^\infty & \text{for some } \bb\in\lt(\tfrac{d+4}{5},+\infty\rt),
		\end{align*}
		and \textbf{one} of the two following conditions
		\begin{align}\label{hyp:grad_3}
			\nabla^2 K &\in L^2 + L^q &\text{ for some } &q\in(1,2),
			\\\label{hyp:besov}
			\nabla K &\in B^1_{1,\infty},
		\end{align}
		and let $f$ be a solution of the \eqref{eq:Vlasov} equation and $\oph$ be a solution of \eqref{eq:Hartree} equation with respective initial conditions 
		\begin{align*}
			f^\init &\in \P\cap L^\infty_{x,\xi} \text{ verifying \eqref{hyp:VP_moment} and \eqref{hyp:VP_intg}}
			\\
			\oph^\init &\in \PP\cap \L^\infty.
		\end{align*}
		Moreover, assume that for a given $n\in 2\N$ such that $n>d$ 
		\begin{equation*}
			\forall \ii\in\Int{1,d},\, \opp_\ii^n\oph^\init\in \L^\infty,
		\end{equation*}
		where $\opp_\ii := -i\hbar\partial_\ii$. Assume also that the initial quantum velocity moment
		\begin{equation}\label{hyp:moments_2}
			M_{n_1}^\init < C \text{ for a given } n_1 \geq \frac{\bb(n-1)+d}{\bb-1},
		\end{equation}
		with $n_1\in 2\N$. Then there exists $T>0$ such that
		\begin{align*}
			M_{n_1} &\in L^\infty([0,T])
			\\
			\opp_\ii^n \oph &\in L^\infty([0,T],\L^\infty) \text{ for any } \ii\in\Int{1,d}
			\\
			\rho_\hbar &\in L^\infty([0,T],L^\infty),
		\end{align*}
		uniformly in $\hbar$, and there exists a constant $C_T$ depending only on the initial conditions and independent of $\hbar$ such that
		\begin{equation*}
			\Wh(f(t),\oph(t)) \leq C_T \lt(\Wh(f^\init,\oph^\init) + \sqrt{\hbar}\rt).
		\end{equation*}
		Moreover, when $\bb \geq \frac{d+2}{3}$ and \eqref{hyp:grad_3} is verified, we can take $T=+\infty$ and the same time estimate as in Theorem~\ref{th:CV} holds. If $\bb \geq \frac{d+2}{3}$ and \eqref{hyp:besov} is verified (which is the case for the Coulomb potential in dimension $d=2$), we obtain the following time dependence instead
		\begin{equation*}
			\Wh(f(t),\oph(t)) \leq \max\lt(\sqrt{d\hbar},\,\Wh(f^\init,\oph^\init)^{e^{t/\sqrt{2}}} e^{\lambda(t)(e^{t/\sqrt{2}}-1)}\rt),
		\end{equation*}
		where
		\begin{align*}
			\lambda(t) &=  C \lt(1 + \|\nabla K\|_{B^1_{1,\infty}}\sup_{[0,t]}(\|\rho\|_{L^\infty}(t) + \|\rho_\hbar\|_{L^\infty}(t))\rt).
		\end{align*}
	\end{thm}
	
	\begin{remark}
		From Theorem \ref{th:comparaison_W2}, we can replace the $\Wh$ pseudo-distance in the left of the semiclassical estimates of the two previous theorems by the classical Wasserstein distance up to adding a constant $\sqrt{2d\hbar}$. Moreover, if the initial states are superposition of coherent states, then we can also replace the $\Wh$ pseudo-distance in the right of the inequalities. This is detailed in Section~\ref{sec:coh_stat}.
	\end{remark}
	
	\begin{remark}
		If $K=\frac{1}{|x|^a}$ or $K = -\ln(|x|)$ if $a=0$ (i.e. $\bb = \frac{d}{a+1}$) and $r=\infty$, we can summarize the results by the following table, where "global" indicates that the result is global in time and "local" that it is proved up to a fixed maximal time. We have highlighted the cases corresponding to the Coulomb interaction.
		\begin{center}
		\begin{tabular}{|rcl|c|c|}
		\hline&&&&\\[-1em]
		&&Settings & Moments & Semiclassical limit \\
		\hline&&&&\\[-1em]
		$d = 2$ & and & $a\in (-1,0]$ & \textbf{global} & \textbf{global} \\
		&&&&\\[-1em]
		$d = 2$ & and & $a\in \lt(0,\frac{1}{2}\rt]$ & global & ? \\[0.1em]
		&&&&\\[-1em]
		$d = 2$ & and & $a\in \lt(\frac{1}{2},\frac{2}{3}\rt]$ & local & ? \\[0.1em]
		\hline&&&&\\[-1em]
		$d = 3$ & and & $a\in \lt(-\frac{1}{2}, \frac{4}{5}\rt)$ & global & global \\[0.2em]
		&&&&\\[-1em] 
		$d = 3$ & and & $a\in\lt[\frac{4}{5},1\rt]$ & \textbf{local} & \textbf{local} \\[0.2em]
		&&&&\\[-1em] 
		$d = 3$ & and & $a\in\lt(1,\frac{8}{7}\rt]$ & local & ? \\[0.2em]
		\hline&&&&\\[-1em]
		$d \geq 4$ & and & $a\in\lt(\frac{d}{2}-2,\frac{2(d-1)}{d+2}\rt)$ & global & global \\[0.4em] 
		&&&&\\[-1em]
		$d \geq 4$ & and & $a\in\lt[\frac{2(d-1)}{d+2},\frac{n(d-1)}{n+d}\rt]$ & local & local \\[0.5em]\hline
		\end{tabular}
		\end{center}
		~\\
		In particular, if $\forall\ii\in\Int{1,d},\opp^4_\ii\oph^\init\in \L^\infty$, it proves the convergence of the Hartree equation with Coulomb interaction potential towards the Vlasov-Poisson equation for short times in dimension $d=3$ and all times for $d=2$ under the assumption that $M_{16}^\init$  is bounded in dimension $d=3$ and that $M_{8}^\init$ is bounded in dimension $d=2$. As an other example, if $a$ is close but smaller than $4/5$ in dimension $d=3$ then \eqref{eq:CV_t} holds as soon as $M_{42}^\init$ is bounded.
	\end{remark}

	\begin{remark}
		The hypothesis $a>\frac{d}{2}-2$ seems harder to remove since it comes from the hypothesis $\nabla^2K\in L^2$ which is needed for the comparison between the negative Sobolev distance and the quadratic Wasserstein distance (see Proposition~\ref{prop:ineq_interpol}).
	\end{remark}
	
	\begin{remark}
		As it can be seen from Proposition~\ref{prop:coupling_estimate_Coulomb}, the semiclassical estimate of Theorem~\ref{th:CV_VP} is actually global in time for the Coulomb potential in dimension $d=3$ provided $\rho\in L^\infty_\loc(\R_+,L^\infty)$. And this would follow from the propagation of order $4$ velocity moments globally in time, since then Theorem~\ref{th:propag_moments} and Proposition~\ref{prop:propag_Lp_m} would imply propagation of higher order moments and weighted Lebesgue norms and the desired bound. In the classical case, global in time propagation of moments is proved in \cite{lions_propagation_1991} through the use of a Duhamel formula in order to use the properties of dispersion of the kinetic transport semigroup. However, we did not manage to use the gain of regularity due to the dispersion. Even if it is possible to express the solution of \eqref{eq:Hartree} through a Duhamel formula for operators, the lack of positivity of the operators involved seems to create difficulties. However, an other effect of dispersion is the decay in time of space moments which we will use in a forthcoming paper to prove global in time estimates for small initial data.
	\end{remark}

	The rest of the paper is organized as follows. In Section~\ref{sec:interp} we generalize the classical kinetic interpolation inequalities which are the key inequalities of our work. In Section~\ref{sec:energy} we recall the conservation of energy and Schatten norms and discuss the case of interaction kernels which do not vanish at infinity.
	
	Sections \ref{sec:moments} and \ref{sec:lebesgue} prove the propagation of quantum moments (Theorem~\ref{th:propag_moments}) and quantum weighted Lebesgue norms uniformly in $\hbar$ (First part of Theorem~\ref{th:CV_VP}). In each case, we first write the classical version of the proof and then the quantum case which is more technical.
	
	In Section~\ref{sec:coupling}, we prove the semiclassical limit in term of the modified Wasserstein distance using the regularity results of previous sections. It finishes the proof of Theorem~\ref{th:CV} and Theorem~\ref{th:CV_VP}.
	
	Finally, Section~\ref{sec:coh_stat} shows that the quantum Lebesgue norms and the quantum Wasserstein pseudo-distance are more natural when looking at superposition of coherent states. It allows us to justify more precisely the definition of the quantum Lebesgue norms and to reformulate our results in terms of the classical Wasserstein distance in this case.

\section{Kinetic quantum interpolation inequalities}\label{sec:interp}

	Let $n\geq0$, $0\leq f=f(x,\xi)\in L^r_{x,\xi}\cap L^1_{x,\xi}(|\xi|^n)$ and $\rho_f =\intd f\d\xi$. Then the classical kinetic interpolation inequality writes
	\begin{align}\label{eq:interp_lp}
		\|\rho_f\|_{L^{p_n}} &\leq C \lt(\intdd f|\xi|^n\d x\d\xi\rt)^{1-\theta} \lt\|f\rt\|_{L^r_{x,\xi}}^{\theta},
	\end{align}
	where $C$ depends only on $d$, $n$ and $r$ and $p'_n = r'+\frac{d}{n}$ and $\theta = \frac{r'}{p'_n}$ with $p'$ denoting the Hölder conjugate of $p$. Even more generally, for $0\leq k\leq n$, we have 
	\begin{align}\label{eq:interp_lp_k}
		\lt\|\intd f|\xi|^k\d \xi\rt\|_{L^{p_{n,k}}} &\leq C \lt(\intdd f|\xi|^n\d x\d\xi\rt)^{1-\theta} \lt\|f\rt\|_{L^r_{x,\xi}}^{\theta},
	\end{align}
	with $p'_{n,k} = r'+\frac{d}{n}$ and $\theta = r'/p'_{n,k}$.
	
	The quantum version of \eqref{eq:interp_lp} is known for $n=2$ and is a variant of Lieb-Thirring inequality (see \cite[(A.6)]{lions_sur_1993}). It reads
	\begin{align}\label{eq:Lieb_Thirring}
		\|\rho\|_{L^p} &\leq C \Tr\lt(-\Delta\op\rt)^{1-\theta} \|\op\|_r^\theta,
	\end{align}
	with $p' = r'+\frac{d}{2}$ and $\theta = \frac{r'}{p'}$.
	It implies \eqref{eq:interp_lp} when $n=2$ by replacing $f$ with $\fh$, even if $\fh$ is not always nonnegative. Recalling the notation $\opp = -i\hbar\nabla$ for the quantum momentum, using the $\L^p$ norm defined by \eqref{def:q_schatten} and remarking that
	\begin{equation*}
		h^{2(1-\theta)}h^{-\theta d/r'} = h^{2-\frac{r'}{r'+d/2}\lt(2+\frac{d}{r'}\rt)} = 1,
	\end{equation*}
	inequality \eqref{eq:Lieb_Thirring} can be written
	\begin{align*}
		\|\rho\|_{L^p} &\leq C \Tr\lt(|\opp|^2\op\rt)^{1-\theta} \|\op\|_{\L^r}^\theta.
	\end{align*}
	By using the results in \cite{egorov_moments_1995}, we obtain the full generalization of \eqref{eq:interp_lp}.
	
	\begin{thm}\label{th:Lieb_Thirring_+}
		Let $n\in 2\N$. Then there exists $C > 0$ depending only on $d$, $r$ and $n$ such that
		\begin{align}\label{eq:Lieb_Thirring_+}
			\|\rho\|_{L^p} &\leq C \Tr\lt(|\opp|^n\op\rt)^{1-\theta} \|\op\|_{\L^r}^\theta,
		\end{align}
		with $p' = r'+\frac{d}{n}$, $\theta = \frac{r'}{p'}$. Moreover, by defining for $k\in 2\N$,
		\begin{align*}
			\rho_{k} &:= \sumj\lambda_j|\opp^\frac{k}{2}\psi_j|^2 = \diag(\opp^\frac{k}{2}\op\cdot\opp^\frac{k}{2}),
		\end{align*}
		for $k<n$, there exists $C > 0$ depending only on $d$, $r$, $n$ and $k$ such that
		\begin{align}\label{eq:Lieb_Thirring_k}
			\|\rho_k\|_{L^\alpha} &\leq C \Tr\lt(|\opp|^n\op\rt)^{1-\theta_k} \|\op\|_{\L^r}^{\theta_k},
		\end{align}
		where $\alpha' = (n/k)'p'$, and $\theta_k = \frac{r'}{\alpha'}$ with $(n/k)'$ denoting the Hölder conjugate of $n/k$.
	\end{thm}
	
	\begin{remark}
		Since for any $u\in\mathcal{D}'(\R^d,\CC)$, $\opp u \in\CC^d$, remark that for any $k\in\N$, $\opp^k u \in \CC^{d^k}$, which leads to $\opp^k u = (\opp_{\ii_1}...\opp_{\ii_k}u)_{(\ii_1,...,\ii_k)\in \Int{1,d}^k}$ and $|\opp^k u|$ is nothing but the natural euclidean norm on $\CC^{d^k}$
		\begin{equation*}
			|\opp^k u|^2 = \sum_{(\ii_1,...,\ii_k)\in \Int{1,d}^k} |\opp_{\ii_1}...\opp_{\ii_d}u|^2.
		\end{equation*}
	\end{remark}
	
	\begin{remark}
		As it can be seen in the proof, when $k=n$, we get an equality in equation~\eqref{eq:Lieb_Thirring_k}
		\begin{equation*}
			\|\rho_n\|_{L^1} = \Tr\lt(|\opp|^n\op\rt) = \intdd \fh |\xi|^n\d x\d\xi.
		\end{equation*}
	\end{remark}

	\begin{remark}
		Taking $\hbar = 1$, we can write \eqref{eq:Lieb_Thirring_k} as
		\begin{align*}
			\|\rho\|_{L^p} &\leq C \Tr\lt((-\Delta)^\frac{n}{2}\op\rt)^{1-\theta} \|\op\|_r^\theta,
		\end{align*}
		which can be written as a \textit{Gagliardo-Nirenberg inequality for orthogonal functions} under the form
		\begin{align*}
			\lt\|\sumj \lambda_j|\psi_j|^2\rt\|_{L^p} &\leq C \lt(\sumj\lambda_j\lt\|\nabla^\frac{n}{2}\psi_j\rt\|_{L^2}^2\rt)^{1-\theta} \lt(\sumj\lambda_j^r\|\psi_j\|_{L^2}^{2r}\rt)^{\theta/r}.
		\end{align*}
	\end{remark}
	
	\begin{demo}[ of Theorem~\ref{th:Lieb_Thirring_+}]
		As proved in \cite[Theorem 1]{egorov_moments_1995}, for any $s>0$ such that $s>1-\frac{d}{n}$, the following bound holds
		\begin{equation}
			\sum_j|\mu_j|^s \leq C_{s,n,d}\intd \mathcal{V}_-^{s+\frac{d}{n}},
		\end{equation}
		where the $\mu_j$ are the negative eigenvalues of $(-\Delta)^\frac{n}{2}+\mathcal{V}$. By taking $\mathcal{V} = -t\rho^{p-1}$ and $s=r'$, the same proof as in \cite{lions_sur_1993} gives inequality \eqref{eq:Lieb_Thirring_+}.
		
		The second inequality requires some more work. We use a vector-valued version of Gagliardo-Nirenberg inequality proved in \cite{schmeiser_vector-valued_2005} which states in particular that for a given Banach space $X$ and any $u\in (H^n\cap L^p)(\R^d,X)$ we have
		\begin{equation}\label{eq:gag_nir_vector}
			\|\nabla^\frac{k}{2} u\|_{L^{2\alpha}(\R^d,X)} \leq C_{d,k,n,p} \|u\|_{L^{2p}(\R^d,X)}^{1-k/n} \|\nabla^\frac{n}{2} u\|_{L^2(\R^d,X)}^{k/n},
		\end{equation}
		for any $(\alpha,p)\in(1,\infty]^2$, $n\in\N$ and $k\leq n$ such that $\frac{1}{\alpha} = \frac{1}{p}\lt(1-\frac{k}{n}\rt)+\frac{k}{n}$. We will use it for the norm given for $\Psi = (\psi_j)_{j\in J}$ by
		\begin{equation*}
			\|\Psi\|_X^2 := \sumj \lambda_j |\psi_j|^2.
		\end{equation*}
		For this norm, by integrating by parts, we remark that
		\begin{align*}
			\|\opp^\frac{n}{2} \Psi\|_{L^2(\R^d,X)}^2 &= \intd \sumj \lambda_j |\opp^{n/2}\psi_j|^2
			\\
			&= \sum_{(j,\ii_1,...,\ii_{n/2})\in J\times\Int{1,d}^{n/2}} \lambda_j \intd \conj{\opp_{\ii_1}...\opp_{\ii_{n/2}}\psi_j}\opp_{\ii_1}...\opp_{\ii_{n/2}}\psi_j
			\\
			&= \sum_{(j,\ii_1,...,\ii_{n/2})\in J\times\Int{1,d}^{n/2}} \lambda_j \intd \conj{\psi_j}\opp_{\ii_1}^2...\opp_{\ii_{n/2}}^2\psi_j
			\\
			&= \sum_{j\in J} \lambda_j \intd \conj{\psi_j}|\opp|^2...|\opp|^2\psi_j
			\\
			&= \Tr\lt(|\opp|^n\op\rt).
		\end{align*}
		Using inequality \eqref{eq:gag_nir_vector} for $\Psi$ and multiplying it by $\hbar^{k/2}$, we obtain
		\begin{align}\label{eq:gag_nir_vector_2}
			\lt\|\rho_{k}\rt\|_{L^\alpha(\R^d,X)} &\leq C_{d,k,n,p} \|\rho\|_{L^p}^{1-\frac{k}{n}} \Tr(|\opp|^n\op)^{k/n},
		\end{align}
		where
		\begin{equation*}
			\frac{1}{\alpha'} = \frac{1}{p'}\lt(1-\frac{k}{n}\rt) = \frac{1}{p'(n/k)'}.
		\end{equation*}	
		Using the first inequality \eqref{eq:Lieb_Thirring_+} to bound $\|\rho\|_{L^p}$ in the left hand side and the fact that $\theta_k = \lt(1-\frac{k}{n}\rt)\theta$, we deduce formula \eqref{eq:Lieb_Thirring_k}.
	\end{demo}
	
\section{Conservation laws}\label{sec:laws}

	In this section, we recall the conservation laws for the \eqref{eq:Vlasov} equation and their equivalent for \eqref{eq:Hartree} equation.
	
\subsection{Conservation of the Schatten norm}

	The Hamiltonian structure of the Vlasov equation implies the preservation of the Lebesgue norms
	\begin{equation*}
		\|f\|_{L^r_{x,\xi}} = \|f^\init\|_{L^r_{x,\xi}}.
	\end{equation*}
	The following property is the quantum equivalent of this conservation law expressed in term of quantum Lebesgue norms.
	
	\begin{prop}\label{prop:propag_Lp}
		Let $\op$ be a solution of the \eqref{eq:Hartree} equation with initial condition $\op^\init \in \PP\cap\L^r$. Then
		\begin{equation*}
			\|\op\|_{\L^r} = \|\op^\init\|_{\L^r}.
		\end{equation*}
	\end{prop}
	
	\begin{demo}
		Assume $r\in\N$. Since $\partial_t\op = [H_\rho,\op]$, we obtain
		\begin{align*}
			\partial_t\op^2 &= \op[H_\rho,\op] + [H_\rho,\op]\op
			\\
			&= [H_\rho,\op^2],
		\end{align*}
		and by an immediate recurrence, for any $n\in\N$, $\partial_t\op^n = [H_\rho,\op^n]$. It implies in particular that
		\begin{align*}
			\dt\Tr(\op^r) &= \Tr([H_\rho,\op^r]) = 0.
		\end{align*}
		Since $\op\geq 0$, we can write $\op = |\op|$ and deduce that $\|\op\|_{\L^r}$ is constant in time. When $r$ is not an integer, the result follows by complex interpolation and the case $r=+\infty$ is obtained by passing to the limit $r\to\infty$.
	\end{demo}
	
\subsection{Conservation of Energy}\label{sec:energy}

	The conservation of energy is a well known property of both \eqref{eq:Vlasov} and \eqref{eq:Hartree} equations, see for example \cite{ginibre_class_1980} and \cite{lions_sur_1993} for the quantum case. For the sake of completeness we write here a short proof with our notations.
	
	\begin{prop}\label{prop:conservation_energy}
		Let $\op\in\PP$ be a solution of \eqref{eq:Hartree} equation. We define the total energy of the system by
		\begin{equation*}
			\Eps_T := M_2 + \intd \rho V,
		\end{equation*}
		where $M_2 = \Tr(|\opp|^2\op)$ and $V = K*\rho$ for a symmetric kernel $K$. Then, as in the classical case, the total energy is conserved
		\begin{equation*}
			\Eps_T(t) = \Eps_T(0).
		\end{equation*}
	\end{prop}
	
	\begin{remark}
		Notice that we can also write
		\begin{equation*}
			\Eps_T = \iintd (|\xi|^2 + V(x)) \fh(x,\xi)\d x\d \xi = \Tr((|\opp|^2 + V)\op),
		\end{equation*}
		which shows that the energy has the same expression with the Wigner transform $\fh$ as in the classical case.
	\end{remark}
	
	\begin{remark}
		By the interpolation inequality \eqref{eq:Lieb_Thirring_+} and assuming that $M_2$ is bounded and $\op\in\L^r\cap\L^1$, we get that $\rho\in L^p$ for $p' \in[r' + d/2,\infty]$. Thus, by Hardy-Littlewood-Sobolev inequality, the negative part of the potential energy $(\Eps_P)_- = \intd \rho V_-$ is bounded for $K_-\in L^{\mathfrak{a},\infty} + L^\infty$ with $\mathfrak{a} = p'/2$. Therefore, if $\op\in\L^r\cap\L^1$ with $r' \leq 2\bb_0 - d/2$,  $(\Eps_P)_-$ is controlled by the kinetic energy and both quantities remains finite if $M_2^\init$ is bounded. It includes the Coulomb interaction in dimension $d=3$. See also \cite{lions_sur_1993}.
		
		If $K$ is not bounded for $|x|\to\infty$ but $K = K_0 + K_\infty \in L^{\mathfrak{a},\infty} + L^\infty(|x|^{-k})$, as in the case of the two-dimensional Coulomb interaction $K = -\ln(|x|)$, $\Eps_P$ can be controlled by assuming for example additional finite space moments 
		\begin{equation*}
			N_k = \intdd \fh(x,\xi) |x|^k\d x\d\xi = \Tr(|x|^k\op) < C.
		\end{equation*}
		In this case, one can indeed write
		\begin{align*}
			\Eps_P = \intd (K_0*\rho) \rho + \intd (K_\infty*\rho) \rho.
		\end{align*}
		The first integral is still controlled as above by $M_2$ if $2\mathfrak{a} \in [r'+d/2,\infty]$. to control the second, we write
		\begin{align*}
			\lt|\intd (K_\infty*\rho) \rho\rt| &\leq C\intd |x-y|^k \rho(\!\d x) \rho(\!\d y)
			\\
			&\leq C\intd (|x|^k+|y|^k) \rho(\!\d x) \rho(\!\d y)
			\\
			&\leq 2C M_0 N_k.
		\end{align*}
		It is easy to see that if $M_2^\init+N_2^\init$ is bounded, then space and velocity moments up to order $2$ remain bounded, since $\partial_t N_2  = \Tr((x\cdot\opp+\opp\cdot x)\op) \leq 2 N_2^{1/2} M_2^{1/2}$, which combined with the conservation of energy leads to
		\begin{align*}
			|\partial_t(M_2 + N_2 + \Eps_P)| &\leq M_2 + N_2
			\\
			&\leq M_2 + N_2 + \Eps_P + C(M_0^{\theta_1}M_2^{\theta_2} + N_0^{\theta_3}N_2^{\theta_4}).
			\\
			&\leq (1+C)(M_2 + N_2 + \Eps_P) + C M_0.
		\end{align*}
		By Gronwall's Lemma and since $\Eps_P$ is controlled by $M_2 + N_2$, we obtain that $M_2 + N_2\in L^\infty_\loc(\R_+)$.
	\end{remark}
	
	\begin{demo}[ of Proposition~\ref{prop:conservation_energy}]
		Since $\Tr(H[H,\op]) = \Tr([H,H]\op) = 0$ and $\partial_t H = \partial_t V$, we obtain
		\begin{align*}
			2\dt\Tr(H\op) &= 2\Tr((\partial_tH) \op) + \frac{2}{i\hbar}\Tr(H [H,\op])
			\\
			&= 2\Tr((\partial_tV)\op)
			\\
			&= 2\intd (\partial_t V) \rho,
		\end{align*}
		and since $K$ is symmetric, we get
		\begin{equation*}
			2\intd (\partial_t V) \rho = \intd (K*\partial_t \rho) \rho + (K*\rho) \partial_t\rho = \dt \intd \rho V.
		\end{equation*}
		Now we remark that
		\begin{equation*}
			2\Tr(H\op) = \Tr(|\opp|^2\op) + 2\Tr(V\op) = M_2 + 2\intd \rho V.
		\end{equation*}
		Thus, we obtain
		\begin{equation*}
			\dt\lt(M_2 + 2\intd \rho V\rt) = 2\dt\Tr(H\op) = \dt \intd \rho V,
		\end{equation*}
		which leads to the result.
	\end{demo}

\section{Propagation of moments}\label{sec:moments}

	We study in this section the propagation independently of $\hbar$ of velocity moments for the Wigner transform of the density operator $\op$ solution of the \eqref{eq:Hartree} equation, which write
	\begin{equation*}
		M_n := \iintd \fh(x,\xi)|\xi|^n\d x\d\xi = \Tr(|\opp|^n\oph).
	\end{equation*}
	To clarify the presentation, we first prove the classical estimate which will be our guideline to prove the semiclassical case.

\subsection{Classical case}

	In this section, we consider only the classical quantities, so that we define
	\begin{align*}
		\rho_n &:= \intd f(x,\xi)|\xi|^n\d\xi
		\\
		M_n &:= \iintd f(x,\xi)|\xi|^n\d x\d\xi = \intd \rho_n.
	\end{align*}
	We can then prove the classical analogue of Theorem~\ref{th:propag_moments}.
	\begin{prop}\label{prop:classique_propag_moments}
		Let $r\geq2$, $\bb_n := \frac{nr'+d}{n+1}$ and $\nabla K\in L^{\bb,\infty}$ for a given $\bb\in[\bb_3,+\infty]$ and $f$ verify the \eqref{eq:Vlasov} equation for $t\in[0,T]$ with initial condition $f^\init\in\P\cap L^r_{x,\xi}$ such that $M_0$ and $M_n^\init$ are bounded for a given $n\geq 2$. Then there exists $T>0$ and $\Phi\in C^0[0,T)$ such that for any $t\in[0,T)$
		\begin{equation}
			M_n \leq \Phi(t).
		\end{equation}
		Moreover, $T=+\infty$ when $\bb \geq \bb_2 = \frac{2r'+d}{3}$.
	\end{prop}

	\begin{demo}[ of Proposition~\ref{prop:classique_propag_moments}]
		Since $M_0$ and $M_n^\init$ are bounded, we deduce that $M_2^\init$ is bounded and by conservation of the energy (Proposition~\ref{prop:conservation_energy}) we deduce that $M_2 \in L^\infty_\loc(\R_+)$. To simplify we write $f = f(t,x,\xi)$. Then we have
		\begin{align*}
			\ddt{M_n} &= \iintd (-\xi\cdot\nabla_x f -E(x)\cdot\nabla_\xi f)|\xi|^n \d x\d\xi
			\\
			&= n\iintd f E(x)\cdot\xi|\xi|^{n-2}\d x\d\xi.
		\end{align*}
		Since $E = -\nabla K * \rho$ with $\nabla K \in L^{\bb,\infty}$, Hölder's and Hardy-Littlewood-Sobolev's inequalities give
		\begin{align}\label{eq:classique_d_t_M_n}
			\lt|\ddt{M_n}\rt| &\leq n\lt\|\intd f |\xi|^{n-1}\d\xi\rt\|_{L^\alpha} \|E\|_{L^{\alpha'}} 
			\\
			&\leq n\|\rho_{n-1}\|_{L^\alpha} \|\rho\|_{L^\beta},
		\end{align}
		where $(\alpha,\beta)\in(1,\infty)^2$ are such that $1+\frac{1}{\alpha'} = \frac{1}{\beta} + \frac{1}{\bb}$ or equivalently
		\begin{equation*}
			\frac{1}{\alpha'} + \frac{1}{\beta'} =  \frac{1}{\bb}.
		\end{equation*}
		By the interpolation inequality \eqref{eq:interp_lp_k}, if we can take $\alpha' = p_{n,n-1}' = np'_n = nr'+d \geq \bb$ and $\theta = \frac{r'}{\alpha'}$, we get
		\begin{equation}\label{eq:classique_rho_n-1_alpha}
			\|\rho_{n-1}\|_{L^\alpha} \leq C M_n^{1-\theta}\|f\|_{L^r_{x,\xi}}^\theta.
		\end{equation}
		Moreover, since the $L^r_{x,\xi}$ is conserved, we can replace $\|f\|_{L^r_{x,\xi}}$ by $\|f^\init\|_{L^r_{x,\xi}}$. If $\beta \leq p_{n-1}$, we can bound $\|\rho\|_{L^\beta}$ using only moments of order less than $n-1$ by using the interpolation inequality \eqref{eq:interp_lp}
		\begin{align*}
			\|\rho\|_{L^\beta} \leq \|\rho\|_{L^{p_{n-1}}}^\frac{p'}{\beta'} \|\rho\|_{L^1}^{1-\frac{p'}{\beta'}}
			\leq C M_0^{1-\frac{p'}{\beta'}}M_{n-1}^{\frac{p'}{\beta'}(1-\frac{r'}{p'_n})} \|f^\init\|_{L^r_{x,\xi}}^\frac{p'r'}{\beta'p'_n}.
		\end{align*}
		Therefore, for $\dt M_n$, the inequality becomes
		\begin{align*}
			\lt|\ddt{M_n}\rt| &\leq C_{d,n,r} M_0^{1-\frac{p'}{\beta'}} \|f^\init\|_{L^r_{x,\xi}}^{\theta + \frac{p'r'}{\beta'p'_n}} M_{n-1}^{\frac{p'}{\beta'}(1-\frac{r'}{p'_n})} M_n^{1-\theta}.
		\end{align*}
		Assuming that $M_{n-1}$ is bounded on $[0,T]$, by Gronwall's Lemma, it implies a bound on $[0,T]$ for $M_n$.
		
		If $\beta > p_{n-1}$, we remark that
		\begin{align*}
			\beta \leq p_n &\ssi \frac{1}{\bb} - \frac{1}{\alpha'} \leq \frac{1}{p'_n}
			\\
			&\ssi
			\frac{1}{\bb}\leq \frac{1}{p'_n}\lt(1+\frac{1}{n}\rt)
			\\
			&\ssi \bb \geq \frac{nr'+d}{n+1} =: \bb_n.
		\end{align*}
		In this case, by interpolation between Lebesgue spaces and by the interpolation inequality \eqref{eq:interp_lp}, we get
		\begin{align*}
			\|\rho\|_{L^\beta} &\leq \|\rho\|_{L^{p_n}}^{\eps} \|\rho\|_{L^{p_{n-1}}}^{1-\eps}
			\\
			&\leq C_{d,n,r} M_{n-1}^{(1-\eps)(1-\theta_{n-1})} M_{n}^{\eps(1-\theta_{n})} \|f^\init\|_{L^r_{x,\xi}}^{(1-\eps)\theta_{n-1}+\eps\theta_{n}},
		\end{align*}
		where $\theta_n = \frac{r'}{p'_n}$ and $\eps\in(0,1)$ is defined by
		\begin{equation}\label{eq:classique_eps0}
			\frac{1}{\beta'} = \frac{\eps}{p'_n} + \frac{1-\eps}{p'_{n-1}}.
		\end{equation}
		By \eqref{eq:classique_d_t_M_n} and \eqref{eq:classique_rho_n-1_alpha}, it implies
		\begin{equation*}
			\lt|\ddt{M_n}\rt| \leq C_{d,n,r} \|f^\init\|_{L^r_{x,\xi}}^{\theta  + (1-\eps)\theta_{n-1} + \eps\theta_{n}} M_{n-1}^{\Theta_0} M_n^{\Theta},
		\end{equation*}
		where
		\begin{align*}
			\Theta_0 &= (1-\eps)(1-\theta_{n-1})
			\\
			\Theta &= 1-\theta + \eps(1-\theta_n).
		\end{align*}
		Using equation \eqref{eq:classique_eps0} to compute $\eps$, we obtain
		\begin{align*}
			\eps &= \lt(\frac{1}{\beta'} - \frac{1}{p'_{n-1}}\rt)\lt(\frac{1}{p'_n} - \frac{1}{p'_{n-1}}\rt)^{-1}
			\\
			&= \lt(\frac{1}{\bb} - \frac{1}{np'_n} - \frac{1}{p'_{n-1}}\rt)\lt(\frac{1}{p'_n} - \frac{1}{p'_{n-1}}\rt)^{-1}
			\\
			&= \frac{n(n-1)}{d} \lt(\frac{p'_np'_{n-1}}{\bb} - \frac{p'_{n-1}}{n} - p'_n\rt)
			\\
			&= \frac{(nr'+d)((n-1)r'+d)}{d\bb} - \frac{(n-1)r'}{d}(1+n) - n.
		\end{align*}
		Since $1-\frac{r'}{p'_n} = \frac{d}{d+nr'}$ we deduce that
		\begin{align*}
			\Theta &= 1-\frac{r'}{np'_n} + \eps \frac{d}{d+nr'}
			\\
			&= \frac{d+(n-1)r'}{d+nr'} + \frac{((n-1)r'+d)}{\bb} - \frac{(n-1)r'}{d+nr'}(1+n) - \frac{nd}{d+nr'}
			\\
			&= \frac{(n-1)r'+d}{\bb} - n + 1 = 1 + n\lt(\frac{\bb_{n-1}}{\bb}-1\rt).
		\end{align*}
		In particular,
		\begin{equation}
			\Theta \leq 1 \ssi \bb \geq \frac{(n-1)r' + d}{n} = \bb_{n-1},
		\end{equation}
		which, by Gronwall's Lemma, allows to prove that $M_n$ is bounded on $[0,T]$ when $M_{n-1}$ is bounded $[0,T]$ and $M_n^\init$ is bounded. In particular, since $M_2$ is bounded by the energy conservation and $\bb = 1 + \frac{d-r'}{n}$ is decreasing with $n$, all the moments will be propagated if the moment of order $3$ is bounded, which is the case when $\bb \geq \frac{d+2r'}{3}$.
	\end{demo}
	
\subsection{Quantum case}

	As we will not consider the \eqref{eq:Vlasov} equation in this section, we will omit to write the $\hbar$ for $\oph$ solution of \eqref{eq:Hartree} equation and $\rho = \diag(\oph)$ to simplify the notations.
	
	\begin{demo}[ of Theorem~\ref{th:propag_moments}]
	To simplify the computations, we define for any $k\in\N$,
		\begin{align*}
			[\opp]^{2k} &:= |\opp|^{2k}
			\\
			[\opp]^{2k+1} &:= |\opp|^{2k}\opp.
		\end{align*}
		\step{1. An inequality for the time derivative of moments} We remark that
		\begin{align*}
			[\opp,H] &= [\opp,V] = -i\hbar\nabla(V\cdot)+i\hbar V\nabla = i\hbar E
			\\
			[|\opp|^2,H] &= \opp\cdot \opp H - H \opp\cdot \opp = \opp\cdot [\opp,H] + [\opp,H]\cdot \opp = i\hbar(\opp\cdot E+E\cdot \opp)
			\\
			[|\opp|^{2n+2},H] &= |\opp|^2\lt[|\opp|^{2n},H\rt] + \lt[|\opp|^{2n},H\rt]|\opp|^2.
		\end{align*}
		By an immediate recurrence, we deduce that for any $n\in\N$,
		\begin{align}\label{eq:pn_H}
			\frac{1}{i\hbar}[|\opp|^{2n+2},H] &= \sum_{k=0}^{n}\binom{n}{k} |\opp|^{2k}(\opp\cdot E+E\cdot \opp)|\opp|^{2(n-k)}.
		\end{align}
		With this formula, we can compute the time derivative of moments as follows
		\begin{align}\nonumber
			\dt\Tr(|\opp|^{2n+2}\op) &= \sum_{k=0}^{n}\binom{n}{k} \Tr\lt(|\opp|^{2k}(\opp\cdot E+E\cdot \opp)|\opp|^{2(n-k)}\op\rt)
			\\\nonumber
			&= \sum_{k=0}^{n}\binom{n}{k} \Tr\lt(([\opp]^{2k+1}\cdot E[\opp]^{2(n-k)}+[\opp]^{2k}E\cdot [\opp]^{2(n-k)+1})\op\rt)
			\\\label{eq:d_t_m_n}
			&= \sum_{k=0}^{2n+1}\binom{n}{\lfloor k/2\rfloor} \Tr\lt([\opp]^{k}\cdot E\cdot [\opp]^{2n+1-k}\op\rt).
		\end{align}
		Recalling that $\op = \sumj \lambda_j\ket{\psi_j}\bra{\psi_j}$, for the term with $k=n$, we have
		\begin{align*}
			 \Tr\lt([\opp]^{n}\cdot E\cdot[\opp]^{n+1}\op\rt)
			 &= \sumj \lambda_j \intd ([\opp]^{n}\conj{\psi_j})\cdot E\cdot([\opp]^{n+1}\psi_j)
			 \\
			 &\leq \intd |E| \rho_{2n}^\frac{1}{2}\rho_{2n+2}^\frac{1}{2}.
		\end{align*}
		For the other terms, for $k<n$, we integrate by parts and use Cauchy-Schwartz inequality to find
		\begin{align}\label{eq:develop_0}
			\Tr\lt([\opp]^{k}\cdot E\cdot[\opp]^{2n+1-k}\op\rt)
			&= \sumj\lambda_j\intd [\opp]^{n-k}(E[\opp]^{k}\conj{\psi_j}) ([\opp]^{n+1}\psi_j)
			\\\nonumber
			&\leq \intd \lt(\sumj\lambda_j\lt|[\opp]^{n-k}(E[\opp]^{k}\psi_j)\rt|^2\rt)^\frac{1}{2} \rho_{2n+2}^\frac{1}{2}.
		\end{align}
		Next we use the definition of $E$ to write
		\begin{align}\label{eq:develop_1}
			\sumj\lambda_j\lt|[\opp]^{n-k}(E[\opp]^{k}\psi_j)\rt|^2 &= \sumj\lambda_j\lt|\sum_{j_2\in J}\lambda_{j_2} [\opp]^{n-k} \lt(\lt(\nabla K*|\psi_{j_2}|^2\rt) [\opp]^{k}\psi_j\rt)\rt|^2.
		\end{align}
		To continue, we introduce the multi-index notation
 		\begin{align*}
			a &= (a_\ii)_{\ii\in\Int{1,d}} \in \N^d \text{ is a finite sequence of integers}
			\\
			|a| &= \sum_{\ii=1}^d a_\ii
			\\
			\opp^a &= (-i\hbar)^{|a|}\partial_{x_1}^{a_1}\partial_{x_2}^{a_2}\dots \partial_{x_d}^{a_d}.
		\end{align*}
		With these notations, we can write
		\begin{align*}
			|\opp|^{2n}(uv) &= \sum_{|a+b|=2n} C_{a,b}^{2n} \, \opp^a(u), \opp^b(v),
		\end{align*}
		where the constants $C_{a,b}^{2n}$ are non-negative integers depending on the multi-indices $a$ and $b$ and such that
		\begin{equation}\label{eq:estim_coefs}
			\sum_{|a+b|=2n} C_{a,b}^{2n} \leq (4d)^n.
		\end{equation}
		More generally, we will write
		\begin{align*}
			[\opp]^n(uv) &= \sum^{\sim}_{|a+b|=n} C_{a,b}^n \, \opp^a(u) \opp^b(v),
		\end{align*}
		where the sum is taken only over the $(a,b)$ such that $|a+b|=n-1$ if $n$ is odd, since then $[\opp]^n(uv)$ is a vector with one free index. Hence, we get
		\begin{align*}
			[\opp]^{n-k} \lt(\lt(\nabla K*|\psi_{j_2}|^2\rt) [\opp]^{k}\psi_j\rt) &= \sum^{\sim}_{|a+b+c|= n-k} C_{a,b,c}^{n-k} \nabla K*\lt(\opp^a(\conj{\psi_{j_2}}) \opp^b(\psi_{j_2})\rt)\opp^c[\opp]^{k}\psi_j.
		\end{align*}
		Moreover, by Cauchy-Schwartz inequality
		\begin{align*}
			\lt|\sum_{j_2\in J}\lambda_{j_2}\opp^a(\conj{\psi_{j_2}}) \opp^b(\psi_{j_2})\rt| &\leq
			\lt(\sumj\lambda_{j}|\opp^a(\psi_{j})|^2\rt)^\frac{1}{2} \lt(\sum_{j_2\in J}\lambda_{j_2}|\opp^b(\psi_{j_2})|^2\rt)^\frac{1}{2}
			\\
			&\leq \rho_{2|a|}^{1/2}\rho_{2|b|}^{1/2}.
		\end{align*}
		Thus, \eqref{eq:develop_1} leads to the following inequality
		\begin{align*}
			\sumj\lambda_j\lt|[\opp]^{n-k}(E[\opp]^{k}\psi_j)\rt|^2
			&\leq
			\sumj\lambda_j\lt|\sum^{\sim}_{|a+b+c|= n-k} C_{a,b,c}^{n-k} \lt( |\nabla K| *\lt(\rho_{2|a|}^{1/2}\rho_{2|b|}^{1/2}\rt) \rt) |\opp^c[\opp]^{k}\psi_j| \rt|^2.
		\end{align*}
		The left hand side can be written under the form
		\begin{equation*}
			\lt\|\sum^{\sim}_{|a+b+c|=n-k}A_{a,b,c}\Psi_c\rt\|_X^2,
		\end{equation*}
		with $\Psi_c = (|\opp^c[\opp]^k\psi_j|)_{j\in J}$ and $\|(u_j)_{j\in J}\|_{X}^2 = \sumj \lambda_j |u_j|^2$. Then, Minkowski's inequality reads
		\begin{equation*}
			\lt\|\sum^{\sim}_{|a+b+c|=n-k}A_{a,b,c}\Psi_c\rt\|_X \leq \sum^{\sim}_{|a+b+c|=n-k}A_{a,b,c}\lt\|\Psi_c\rt\|_X.
		\end{equation*}
		Remarking that $\|\Psi_c\|_X^2 \leq \rho_{2|c|+2k}$, we obtain
		\begin{align*}
			\sumj\lambda_j\lt|[\opp]^{n-k}(E[\opp]^{k}\psi_j)\rt|^2
			&\leq \lt(\sum^{\sim}_{|a+b+c|= n-k} C_{a,b,c}^{n-k} \lt(|\nabla K| *\lt(\rho_{2|a|}^{1/2}\rho_{2|b|}^{1/2}\rt) \rt) \rho_{2|c|+2k}^{1/2}\rt)^2.
		\end{align*}
		Combining this inequality with \eqref{eq:develop_0} and \eqref{eq:estim_coefs}, we obtain
		\begin{align*}
			\Tr\lt([\opp]^k\cdot E\cdot [\opp]^{2n+1-k}\op\rt) &\leq \intd \sum^{\sim}_{|a+b+c|= n-k} C_{a,b,c}^{n-k} \lt(\nabla K*\lt(\rho_{2|a|}^{1/2}\rho_{2|b|}^{1/2}\rt) \rt) \rho_{2|c|+2k}^\frac{1}{2} \rho_{2n+2}^\frac{1}{2}
			\\
			&\leq (4d)^\frac{n-k}{2}\sup_{|a+b+c|= n} \lt\|\nabla K*\lt(\rho_{2|a|}^{1/2}\rho_{2|b|}^{1/2}\rt) \rho_{2|c|}^\frac{1}{2}\rt\|_{L^2}
			\|\rho_{2n+2}\|_{L^1}^\frac{1}{2}
			\\
			&\leq (4d)^\frac{n-k}{2}C_K\sup_{|a+b+c|= n} \lt\|\rho_{2|a|}\rt\|_{L^\alpha}^\frac{1}{2} \lt\|\rho_{2|b|}\rt\|_{L^\beta}^\frac{1}{2} \lt\|\rho_{2|c|}\rt\|_{L^\gamma}^\frac{1}{2} M_{2n+2}^\frac{1}{2},
		\end{align*}
		where $C_K = \|\nabla K\|_{L^{\bb,\infty}}$,
		\begin{equation}\label{eq:b}
			\frac{1}{\alpha'}+\frac{1}{\beta'}+\frac{1}{\gamma'} = \frac{2}{\bb},
		\end{equation}
		and we used Hölder's inequality and the weak Young's inequality. The case $k>n$ is treated in the same way. Thus, from \eqref{eq:d_t_m_n} and the identity
		\begin{align*}
			\sum_{k=0}^{2n+1}\binom{n}{\lfloor k/2\rfloor} \lt(4d\rt)^\frac{n-k}{2} &= \lt(2\sqrt{d}\rt)^n\lt(\sum_{k=0}^{n}\binom{n}{k} \lt(2\sqrt{d}\rt)^{-(2k+1)} + \sum_{k=0}^{n}\binom{n}{k} \lt(2\sqrt{d}\rt)^{-2k}\rt)
			\\
			&= \frac{\lt(1+2\sqrt{d}\rt)(1+4d)^n}{(4d)^\frac{n+1}{2}} =: c_{d,2n+2},
		\end{align*}
		we deduce that for any $n\in\N$,
		\begin{align*}
			\dt\Tr(|\opp|^{2n+2}\op) &\leq c_{d,2n+2}C_K \sup_{|a+b+c|= n} \lt\|\rho_{2|a|}\rt\|_{L^\alpha}^\frac{1}{2} \lt\|\rho_{2|b|}\rt\|_{L^\beta}^\frac{1}{2} \lt\|\rho_{2|c|}\rt\|_{L^\gamma}^\frac{1}{2} M_{2n+2}^\frac{1}{2}.
		\end{align*}
		\step{2. Using the kinetic interpolation} To simplify the notations, we will fix $n\in2\N$ and write previous formula as
		\begin{equation}\label{eq:d_t_m_n_2}
			\ddt{M_n} \leq c_{d,n}C_K M_{n}^\frac{1}{2} \sup_{|a+b+c|= n/2-1} \lt\|\rho_{2|a|}\rt\|_{L^\alpha}^\frac{1}{2} \lt\|\rho_{2|b|}\rt\|_{L^\beta}^\frac{1}{2} \lt\|\rho_{2|c|}\rt\|_{L^\gamma}^\frac{1}{2}.
		\end{equation}
		To bound the right term by powers of $M_n$, we use the kinetic quantum interpolation inequalities \eqref{eq:Lieb_Thirring_k}, which gives for any $k \in \{2|a|,2|b|,2|c|\} \subset \Int{0,n-2}$
		\begin{align}
			\|\rho_k\|_{L^{p_n(k)}} &\leq
			C_{d,r,n,k} M_n^{1-\theta_n(k)} \|\op\|_{\L^r}^{\theta_n(k)},
		\end{align}
		where $p_{n}'(k) = (n/k)'p'_n$ with $p'_n = r'+\frac{d}{n}$ and $\theta_n(k) = \frac{r'}{p'_n(k)}$. Since $k\leq n-2$ the same inequality holds by replacing $n$ by $n-2$. If we can choose $\alpha,\beta,\gamma>1$ and $\eps\in(0,1)$ such that
		\begin{align}\label{eq:interp_a}
			\frac{1}{\alpha'} &= \frac{\eps}{p_{n}'(2|a|)} + \frac{1-\eps}{p_{n-2}'(2|a|)}
			\\\label{eq:interp_b}
			\frac{1}{\beta'} &= \frac{\eps}{p_{n}'(2|b|)} + \frac{1-\eps}{p_{n-2}'(2|b|)}
			\\\label{eq:interp_c}
			\frac{1}{\gamma'} &= \frac{\eps}{p_{n}'(2|c|)} + \frac{1-\eps}{p_{n-2}'(2|c|)}.
		\end{align}
		By interpolation and since by Proposition~\ref{prop:propag_Lp}, $\|\op\|_{\L^r} = \|\op^\init\|_{\L^r}$, we get 
		\begin{align*}
			\|\rho_{2|a|}\|_{L^\alpha} &\leq \|\rho_{2|a|}\|_{L^{p_{n}(2|a|)}}^\eps \|\rho_{2|a|}\|_{L^{p_{n-2}(2|a|)}}^{1-\eps}
			\\
			&\leq C_{d,r,n,|a|} \|\op^\init\|_{\L^r}^{C_{a,n,\eps}} M_{n-2}^{(1-\eps)(1-\theta_{n-2}(2|a|))} M_{n}^{\eps(1-\theta_{n}(2|a|))}.
		\end{align*}
		Since $|a+b+c|=n/2-1$, we remark that
		\begin{align*}
			\frac{1}{2}\lt(\frac{1}{p_{n-2}'(2|a|)} + \frac{1}{p_{n-2}'(2|b|)} + \frac{1}{p_{n-2}'(2|c|)}\rt)
			&= \frac{1}{2p_{n-2}'}\lt(3-2\frac{|a|+|b|+|c|}{n-2}\rt)
			\\
			&= \frac{1}{p_{n-2}'} = \frac{n-2}{(n-2)r'+d}
			\\
			\frac{1}{\bb_{n}} := \frac{1}{2}\lt(\frac{1}{p_{n}'(2|a|)} + \frac{1}{p_{n}'(2|b|)} + \frac{1}{p_{n}'(2|c|)}\rt)
			&= \frac{1}{2p_{n}'}\lt(3-2\frac{|a|+|b|+|c|}{n}\rt)
			\\
			&= \frac{1}{p_{n}'}\lt(1+\frac{1}{n}\rt) = \frac{n+1}{nr'+d}.
		\end{align*}
		Therefore, by \eqref{eq:b}, we get
		\begin{equation}\label{eq:epsilon}
			\frac{1}{\bb} = \frac{\eps}{\bb_n} + \frac{1-\eps}{p'_{n-2}}.
		\end{equation}
		Let first assume that $\bb \leq p'_{n-2}$. Then, since by assumption $\bb \geq \bb_n$, we can find $(\alpha,\beta,\gamma,\eps)$ verifying \eqref{eq:interp_a}, \eqref{eq:interp_b} and \eqref{eq:interp_c}. Hence, \eqref{eq:d_t_m_n_2} becomes
		\begin{align}
			\ddt{M_n} &\leq C_{d,r,n}C_K \|\op^\init\|_{\L^r}^{\Theta_2} M_{n-2}^{\Theta_0} M_{n}^{\frac{1}{2}+\Theta_1},
		\end{align}
		with
		\begin{align*}
			\Theta_1 &= \frac{\eps}{2}\lt(3 - \theta_{n}(a) - \theta_{n}(b) - \theta_{n}(c)\rt)
			= \eps\lt(\frac{3}{2} - \frac{r'}{\bb_{n}}\rt)
			\\
			\Theta_0 &= (1-\eps)\lt(\frac{3}{2} - \frac{r'}{p'_{n-2}}\rt)
			\\
			\Theta_2 &= \frac{3}{2} - \Theta_1 - \Theta_0.
		\end{align*}
		From \eqref{eq:epsilon}, we can compute $\eps$ and we get
		\begin{equation*}
			\eps = \frac{nr'+d}{(n-2)r'+3d}\lt(\frac{(n-2)r'+d}{\bb}-(n-2)\rt).
		\end{equation*}
		It leads to the following formula for $\Theta = 1/2 + \Theta_1$
		\begin{equation*}
			\Theta = 1 + \frac{d+(n-2)r'}{2}\lt(\frac{1}{\bb} - \frac{n-1}{d+(n-2)r'}\rt) = 1 + \frac{n-1}{2}\lt(\frac{\bb_{n-2}}{\bb} - 1\rt).
		\end{equation*}
		In particular,
		\begin{align*}
			\Theta &\leq 1 \ \ssi \bb_{n-2} \leq \bb
			\\
			\Theta &\underset{n\to\infty}{\rightarrow} 1.
		\end{align*}
		The result then follows by Gronwall's Lemma. If $\bb > p'_{n-2}$, it is no more possible to write \eqref{eq:epsilon}, but we can still find $\tilde{\eps}\in(0,1)$ such that
		\begin{equation*}
			\frac{1}{\bb} = \frac{1-\tilde{\eps}}{\tilde{\bb}} + \tilde{\eps}\lt(\frac{\eps}{\bb_n} + \frac{1-\eps}{p'_{n-2}}\rt),
		\end{equation*}
		where
		\begin{equation*}
			\frac{1}{\tilde{\bb}} = \frac{1}{p'_{2|a|}(2|a|)} + \frac{1}{p'_{2|b|}(2|b|)} + \frac{1}{p'_{2|c|}(2|c|)} = 0,
		\end{equation*}
		and we obtain
		\begin{align}
			\ddt{M_n} &\leq C_n(M_{2|a|},M_{2|b|},M_{2|c|}) M_{n-2}^{\tilde{\eps}\Theta_0} M_n^{\frac{1}{2}+\tilde{\eps}\Theta_1},
		\end{align}
		with $\tilde{\Theta} = 1/2 + \tilde{\eps} \Theta_1 \leq \Theta$ and we can again conclude by Gronwall's Lemma.
	\end{demo}
	
\section{Propagation of higher Lebesgue weighted norms}\label{sec:lebesgue}

\subsection{Classical case}

	As previously, we first do the proof in the classical case as a guideline for the proof of the quantum case. The goal here is to propagate $\|f\|_{L^p_{x,\xi}(|\xi|^n)}$ norms uniformly in $p$. Together with the uniform bound on $\|f\|_{L^p_{x,\xi}}$, it leads to the following bound for some $C,T>0$ and any $t\in[0,T]$
	\begin{equation*}
		0 \leq f(t,x,\xi) \leq \frac{C}{1+|\xi|^n}.
	\end{equation*}
	For $n>d$, this bound implies that $\rho := \intd f\d\xi \in L^\infty([0,T]\times\R^d)$.
	
	\begin{prop}\label{prop:classique_propag_Lp_m}
		Assume $E\in L^\infty([0,T],L^\infty)$ and let $f$ be a solution of the \eqref{eq:Vlasov} equation such that $f^\init\in L^p$ and $f^\init|\xi|^n\in L^p$ for a given $p\in [1,\infty]$. Then
		\begin{equation*}
			\lt\|f|\xi|^n\rt\|_{L^p_{x,\xi}} \leq \lt(\lt\|f^\init|\xi|^n\rt\|_{L^p_{x,\xi}}^\frac{1}{n} + \|E\|_{L^\infty}\lt\|f^\init\rt\|_{L^p_{x,\xi}}^\frac{1}{n} t\rt)^n.
		\end{equation*}
	\end{prop}
	
	\begin{cor}
		Assume $f$ verifies the hypothesis of Proposition~\ref{prop:classique_propag_Lp_m} for $n>d$ and $p=\infty$. Then $\rho\in L^\infty((0,T),L^\infty)$ and
		\begin{align*}
			\|\rho\|_{L^\infty} &\leq C \|f\|_{L^\infty_{x,\xi}(1+|\xi|^n)}
			\\
			&\leq \lt(\lt\|f^\init|\xi|^n\rt\|_{L^p_{x,\xi}}^\frac{1}{n} + \|E\|_{L^\infty}\lt\|f^\init\rt\|_{L^\infty_{x,\xi}}^\frac{1}{n} t\rt)^n + \|f^\init\|_{L^\infty}.
		\end{align*}
	\end{cor}
	
	\begin{demo}
		Since $f = f(t,x,v)$ is solution of the \eqref{eq:Vlasov} equation, differentiating with respect to time and integrating by parts, we get
		\begin{align*}
			\frac{1}{p}\dt \iintd \lt|f|\xi|^n\rt|^p \d x\d \xi &= \iintd \lt|f\rt|^{p-2}f(-\xi\cdot\nabla_x f - E(x)\cdot\nabla_\xi f) |\xi|^{np}\d x\d\xi
			\\
			&= \frac{1}{p}\iintd (-\xi\cdot\nabla_x \lt(|f|^p\rt) - E(x)\cdot\nabla_\xi \lt(|f|^p\rt)) |\xi|^{np}\d x\d\xi
			\\
			&= n\iintd |f|^p E(x)\cdot\xi|\xi|^{np-2}\d x\d\xi.
		\end{align*}
		Using the fact that $E\in L^\infty$ and Hölder's inequality, we obtain
		\begin{align*}
			\frac{1}{p}\dt \lt\|f|\xi|^n\rt\|_{L^p_{x,\xi}}^p &\leq n \|E\|_{L^\infty} \iintd |f|^p |\xi|^{np-1}\d x\d\xi
			\\
			&\leq n \|E\|_{L^\infty} \lt\|f|\xi|^n\rt\|_{L^p_{x,\xi}}^{p-\frac{1}{n}} \lt\|f\rt\|_{L^p_{x,\xi}}^\frac{1}{n}.
		\end{align*}
		This inequality can be written
		\begin{align*}
			\dt\lt\|f|\xi|^n\rt\|_{L^p_{x,\xi}} &\leq n \|E\|_{L^\infty} \lt\|f|\xi|^n\rt\|_{L^p_{x,\xi}}^{1-\frac{1}{n}} \lt\|f\rt\|_{L^p_{x,\xi}}^\frac{1}{n}.
		\end{align*}
		Then by conservation of the $L^p_{x,\xi}$ norm and Gronwall's Lemma, we deduce that
		\begin{align*}
			\lt\|f|\xi|^n\rt\|_{L^p_{x,\xi}} \leq \lt(\lt\|f^\init|\xi|^n\rt\|_{L^p_{x,\xi}}^\frac{1}{n} + \|E\|_{L^\infty}\lt\|f^\init\rt\|_{L^p_{x,\xi}}^\frac{1}{n} t\rt)^n,
		\end{align*}
		and if $\lt\|f\rt\|_{L^\infty_{x,\xi}} < \infty$ and $\lt\|f^\init|\xi|^n\rt\|_{L^\infty_{x,\xi}} < \infty$, we can pass to the limit $p\to\infty$.
	\end{demo}
	
\subsection{Quantum case}

	In this section, we again only focus on the quantum objects, so that we will write $\op := \oph$ and $\rho := \diag(\op)$ to simplify the notations. For $k\in\R_+$, we define the $\L(|\opp|^k)$ space as the space of compact operators $\op$ such that
	\begin{equation*}
		\|\op\|_{\L^p(|\opp|^k)} := \||\opp|^k\op\|_{\L^p} < C,
	\end{equation*}
	where $\L^p$ is defined by \eqref{def:q_schatten}. Remark that if $\op$ is self-adjoint, then $|\op|\opp|^k|^2 = |\opp|^k|\op|^2|\opp|^k$ and by cyclicity of the trace, for any $p\in2\N$, 
	\begin{equation*}
		\|\op\|_{\L^p(|\opp|^k)} = \||\op|\opp|^k\|_{\L^p}.
	\end{equation*}
	Actually, as proved in \cite{erdos_derivation_2001}, this is true also for $p=1$ and can be easily generalized to any $p\in \R_+$, since for any self-adjoint compact operators $A$ and $B$, as pointed out in \cite[Formula~(1.3)]{simon_trace_2005}, the singular values are the same for $AB$ and $(AB)^* = BA$, which leads to
	\begin{equation}\label{eq:commut}
		\|AB\|_p = \|BA\|_p.
	\end{equation}
	We recall Hölder's inequality (see e.g. \cite[Theorem~2.8]{simon_trace_2005}) which reads for any compact operators $A$ and $B$
	\begin{align}\label{eq:Holder_op}
		\|AB\|_r \leq \|A\|_p \|B\|_q \text{ when } \frac{1}{p} + \frac{1}{q} = \frac{1}{r},
	\end{align}
	and the Araki-Lieb-Thirring inequality \cite[Theorem~1]{araki_inequality_1990} which reads
	\begin{align}\label{eq:ALT}
		\Tr\lt((BAB)^{qr}\rt) \leq \Tr\lt((B^qA^qB^q)^r\rt),
	\end{align}
	for any operators $A,B\geq 0$ and $(q,r)\in[1,\infty)\times\R_+$. Remark that for $A,B\geq 0$, since $|AB| = (BA^2B)^\frac{1}{2}$, we can rewrite \eqref{eq:ALT} as
	\begin{equation}\label{eq:ALT_3}
		\|AB\|_{qr}^q \leq \|A^qB^q\|_r \text{ for any } q\geq 1.
	\end{equation}
	
	From these inequalities we deduce the following interpolation inequality
	\begin{prop}
		Let $A\geq 0$ be a compact operator, then for any $\theta\in[0,1]$
		\begin{equation*}
		\|AB^\theta\|_p \leq \|AB\|_p^\theta \|A\|_p^{1-\theta}.
		\end{equation*}
	\end{prop}
	
	\begin{demo}
		Since $A\geq 0$, we can write $A = A^\theta A^{1-\theta}$ and by Hölder's inequality \eqref{eq:Holder_op}, we obtain
		\begin{align*}
			\lt\|AB^\theta\rt\|_p &\leq
			\lt\|A^\theta B^\theta\rt\|_{p/\theta}
			\lt\|A^{1-\theta}\rt\|_{p/(1-\theta)}
			\\
			&\leq \lt\|A^\theta B^\theta\rt\|_{p/\theta} \lt\|A\rt\|_p^{1-\theta}.
		\end{align*}
		Then, we use \eqref{eq:ALT_3} with $q=1/\theta\geq 1$ to get
		\begin{align*}
			\lt\|B^\theta A^\theta\rt\|_{p/\theta} &\leq \lt\|AB\rt\|_p^\theta,
		\end{align*}
		which proves the result.
	\end{demo}
	
	As a corollary of the previous proposition, taking $B=|\opp|^n$ and $A=\op$, we obtain results for the $\L(|\opp|^k)$ norm.
	
	\begin{cor}
		Let $\op$ be a nonnegative hermitian operator, then for any $0\leq k\leq n<\infty$
		\begin{align}\label{eq:interp_lp_m}
			\|\op\|_{\L^p(|\opp|^k)} &\leq \|\op\|_{\L^p(|\opp|^n)}^{k/n} \|\op\|_{\L^p}^{1-k/n}.
		\end{align}
	\end{cor}
	
	We are now ready to prove the propagation of weighted quantum Schatten norms.
	
	\begin{prop}\label{prop:propag_Lp_m}
		Let $\opp_\ii := -i\hbar\partial_\ii$ for a given $\ii\in\Int{1,d}$,
		\begin{align*}
			\nabla K &\in L^\bb + L^\infty & \text{for some } \bb\in(1,+\infty),
		\end{align*}
		$r\in(\bb',\infty]$ and $\op\in \PP\cap\L^r$ verify \eqref{eq:Hartree} equation. Assume moreover that $M_{n_1}$ is bounded on $[0,T]$ for a given $T>0$ and a given $n_1\in \N$ and that $\op^\init\in \L^{2p}(\opp_\ii^n)$ for a given $p\in\N\cup\{\infty\}$ such that $2p\leq r$ and a given $n\in \N$ such that
		\begin{equation}\label{eq:hyp_n}
			n \leq \theta\,n_1  + 1-\frac{d}{\bb},
		\end{equation}
		with $\theta = 1-\frac{r'}{\bb}$. Then for any $t\in[0,T]$,
		\begin{equation}\label{eq:propag_Lp_m}
			\|\op\|_{\L^{2p}(\opp_\ii^n)} \leq 2^n \lt(\|\op^\init\|_{\L^{2p}(\opp_\ii^n)} + \tilde{C}_{\op^\init} \lt(t+\int_0^tM_{n_1}^\theta\rt)^n\rt),
		\end{equation}
		where $\tilde{C}_{\op^\init} = (4^n C_{d,r,n_1} \|\nabla K\|_{L^\bb} (1+M_0))^n \lt\|\op^\init\rt\|_{\L^{2p}}^{1+\frac{nr'}{\bb}}$. In particular, for $r=p=\infty$, we obtain
		\begin{equation}\label{eq:propag_Linfty_m}
			\|\op\|_{\L^\infty(\opp_\ii^n)} \leq 2^n \lt(\|\op^\init\|_{\L^\infty(\opp_\ii^n)} + \tilde{C}_{\op^\init} \lt(t+\int_0^tM_{n_1}^\theta\rt)^n\rt),
		\end{equation}
		with $\tilde{C}_{\op^\init} = (4^n C_{d,n_1} \|\nabla K\|_{L^\bb} (1+M_0))^n \lt\|\op^\init\rt\|_{\L^\infty}^{1+\frac{n}{\bb}}$.
	\end{prop}
	
	\begin{cor}\label{cor:rho_L_infty}
		With the hypotheses of Proposition~\ref{prop:propag_Lp_m}, assume that for a given $n>d/2$, $\op^\init\in \L^\infty(\opp_\ii^{2n})$ for all $\ii\in\Int{1,d}$. Then
		\begin{equation}\label{eq:rho_L_infty}
			\|\rho\|_{L^\infty} \leq c_{d,n} \|\op\|_{L^\infty((0,T),\L^{\infty}(1+\opp^{2n}))},
		\end{equation}
		which is bounded independently from $\hbar$.
	\end{cor}
	
	\begin{demo}[ of Corollary~\ref{cor:rho_L_infty}]
		To prove \eqref{eq:rho_L_infty}, we remark that from Proposition~\ref{prop:propag_Lp_m},
		\begin{equation*}
			P_n\op := \lt(1+\sum_{\ii=1}^{d} \opp_\ii^{2n}\rt)\op \in L^\infty([0,T],\L^\infty).
		\end{equation*}
		Since $P_n$ and $\rho$ are nonnegative self-adjoint operators, by using \eqref{eq:ALT_3} for $r=\infty$ and $q=2$, we obtain
		\begin{equation*}
			\|P_n^\frac{1}{2}\op P_n^\frac{1}{2}\|_\infty = \|\op^\frac{1}{2} P_n^\frac{1}{2}\|^2_\infty \leq \|P_n\op\|_{\infty} \leq C_\rho h^d,
		\end{equation*}
		where $C_\rho = \|P_n\op\|_{L^\infty((0,T),\L^\infty)}$. From this, we get that for any $\varphi\in L^2$,
		\begin{equation*}
			\langle \varphi|P_n^\frac{1}{2}\op P_n^\frac{1}{2}\varphi\rangle \leq C_\rho h^d\|\varphi\|^2 = \langle \varphi | C_\rho h^d\varphi\rangle,
		\end{equation*}
		or equivalently $P_n^\frac{1}{2}\op P_n^\frac{1}{2} \leq C_\rho h^d$. It implies that $A := C_\rho h^d - P_n^\frac{1}{2}\op P_n^\frac{1}{2}$ is a nonnegative self-adjoint operator. Using the Fourier transform, we remark that $P_n$ is invertible and that for any $\varphi\in L^2$ we have
		\begin{align*}
			P_n^{-1}\varphi(x) &= \F_y\lt(\frac{\hat{\varphi}}{1+\sum_{\ii=1}^{d}|h y_\ii|^{2n}}\rt)(x)
			\\
			&= \intd\F_y\lt(\frac{1}{1+\sum_{\ii=1}^{d}|h y_\ii|^{2n}}\rt)(x-z) \varphi(z)\d z.
		\end{align*}
		Since $P_n^{-\frac{1}{2}}$ is a positive operator, we deduce that $C_\rho h^d P_n^{-1} - \op = P_n^{-\frac{1}{2}} A P_n^{-\frac{1}{2}}$ is a nonnegative operator of diagonal
		\begin{align*}
			0 \leq k(x,x)& = C_\rho h^d\F_y\lt(\frac{1}{1+\sum_{\ii=1}^{d}|h y_\ii|^{2n}}\rt)(0) - \rho(x)
			\\
			&= C_\rho h^d \intd \frac{1}{1+\sum_{\ii=1}^{d}|h y_\ii|^{2n}}\d y - \rho(x)
			\\
			&= c_{d,n} C_\rho  - \rho(x),
		\end{align*}
		where, since $2n>d$,
		\begin{equation*}
			c_{d,n} := \intd \frac{\d x}{1+\sum_{\ii=1}^{d}| x_\ii|^{2n}} < \infty.
		\end{equation*} 
		Since $\rho\geq 0$, we deduce that
		\begin{equation*}
			0 \leq \rho(x) \leq c_{d,n}C_\rho,
		\end{equation*}
		which proves the result.
	\end{demo}
	
	\begin{demo}[ of Proposition~\ref{prop:propag_Lp_m}]
		By cyclicity of the trace, for $p\in\N$, $\ii\in\Int{1,d}$ and $n\in\N$ we have
		\begin{align*}
			\Tr(\lt| \opp_\ii^n \op \rt|^{2p}) &= \Tr\lt(\lt( \opp_\ii^{2n} \op^2\rt)^{p}\rt).
		\end{align*}
		Therefore, using again the cyclicity of the trace
		\begin{align*}
			\frac{i\hbar}{p}\dt \Tr(\lt| \opp_\ii^n \op \rt|^{2p}) &=  \Tr\lt(\lt( \opp_\ii^{2n} \op^2\rt)^{p-1}\opp_\ii^{2n}[H,\op^2]\rt)
			\\
			&= \Tr\lt(\lt( \opp_\ii^{2n} \op^2\rt)^{p-1}\opp_\ii^{2n}H\op^2\rt) - \Tr\lt(\lt( \opp_\ii^{2n} \op^2\rt)^{p-1}\opp_\ii^{2n}\op^2 H\rt)
			\\
			&= \Tr\lt(\opp_\ii^{2n}H\op^2 \lt( \opp_\ii^{2n} \op^2\rt)^{p-1}\rt) - \Tr\lt( H\opp_\ii^{2n}\op^2 \lt(\opp_\ii^{2n} \op^2\rt)^{p-1} \rt)
			\\
			&= \Tr\lt([\opp_\ii^{2n},H]\op^2 \lt( \opp_\ii^{2n} \op^2\rt)^{p-1}\rt)
			\\
			&= \Tr\lt(\op[\opp_\ii^{2n},H]\op \lt| \opp_\ii^n \op\rt|^{2p-2}\rt).
		\end{align*}
		Now we write $[\opp_\ii^{2n},H]$ in terms of $E$ thanks to formula \eqref{eq:pn_H} to obtain in the same way
		\begin{align}\nonumber
			\frac{1}{p}\dt \Tr(\lt| \opp_\ii^n \op \rt|^{2p}) &= \sum_{k=0}^{n-1}\binom{n-1}{k} \Tr\lt(\op\opp_\ii^{2(n-1-k)}(\opp_\ii E_\ii+E_\ii \opp_\ii)\opp_\ii^{2k}\op P\rt)
			\\\label{eq:d_t_lp_m}
			&= \sum_{k=1}^{2n}\binom{n-1}{\lfloor (k-1)/2\rfloor} \Tr\lt(\op\opp_\ii^{2n-k} E_\ii \opp_\ii^{k-1}\op P\rt),
		\end{align}
		where $P = \lt| \opp_\ii^n \op\rt|^{2p-2}$.
		Since $\Tr(A^*) = \conj{\Tr(A)}$, the following holds
		\begin{align*}
			\Tr\lt(\op\opp_\ii^{2n-k} E_\ii \opp_\ii^{k-1}\op P\rt) = \conj{\Tr\lt(\op\opp_\ii^{k-1} E_\ii \opp_\ii^{2n-k}\op P\rt)},
		\end{align*}
		so that \eqref{eq:d_t_lp_m} becomes
		\begin{align}\label{eq:d_t_lp_m_2}
			\frac{1}{p}\dt \Tr(\lt| \opp_\ii^n \op \rt|^{2p}) &= 2\Re\sum_{k=1}^{n}\binom{n-1}{\lfloor (k-1)/2\rfloor} \Tr\lt(\op\opp_\ii^{2n-k} E_\ii \opp_\ii^{k-1}\op P\rt).
		\end{align}
		To treat the right term, we remark that Leibniz rule for differentiation leads to
		\begin{align*}
			\opp_\ii^{n-k}E_\ii = \sum_{m=0}^{n-k} \binom{n-k}{m} (\opp_\ii^{m}(E_\ii))\opp_\ii^{n-k-m}.
		\end{align*}
		Therefore we obtain
		\begin{align}\label{eq:util_leibniz}
			\Tr\lt(\op\opp_\ii^{2n-k} E_\ii \opp_\ii^{k-1}\op P\rt) &= \sum_{m=0}^{n-k}\binom{n-k}{m} \Tr\lt(\op\opp_\ii^n (\opp_\ii^{m}(E_\ii)) \opp_\ii^{n-m-1}\op P\rt).
		\end{align}
		Thus we can use Hölder's inequality \eqref{eq:Holder_op} and the interpolation inequality \eqref{eq:interp_lp_m} to get
		\begin{align}\nonumber
			|\Tr(\op\opp_\ii^n \lt(\opp_\ii^{m}(E_\ii)\rt) &\opp_\ii^{n-m-1}\op P)|
			\\\nonumber
			&=
			\lt|\Tr\lt((\opp_\ii^{m}(E_\ii)) \opp_\ii^{n-m-1}\op P\op\opp_\ii^n\rt)\rt|
			\\\nonumber
			&\leq \|\opp_\ii^{m}(E_\ii)\|_\infty \lt\|\opp_\ii^{n-m-1}\op\rt\|_{2p} \lt\||\opp_\ii^n\op|^{2p-2}\rt\|_{\frac{2p}{(2p-2)}} \lt\|\opp_\ii^{n}\op\rt\|_{2p}
			\\\label{eq:util_holder}
			&\leq \|\opp_\ii^{m}(E_\ii)\|_{L^\infty}  \lt\|\op\rt\|_{2p}^\frac{m+1}{n} \lt\|\opp_\ii^{n}\op\rt\|_{2p}^{2p-\frac{m+1}{n}}.
		\end{align}
		The term $\|\op\|_{2p}$ will be controlled by propagation of the $\L^p$ norm (see Proposition~\ref{prop:propag_Lp}). To control $\|\opp_\ii^{m}(E_\ii)\|_{L^\infty}$ for any $m\in\Int{0,n-1}$, by interpolation, it is sufficient to prove that it is bounded for $m=0$ and $m=n-1$. We use again the Leibniz rule to get
		\begin{align*}
			-\opp_\ii^{m}(E_\ii) &= \opp_\ii^{m}\lt(\nabla K * \rho\rt)
			\\
			&= \nabla K * \sum_j \lambda_j \opp_\ii^{m}(|\psi_j|^2)
			\\
			&= \sum_{l=0}^m \binom{m}{l}\nabla K * \sum_j \lambda_j \opp_\ii^{l}(\conj{\psi_j})\opp_\ii^{m-l}(\psi_j).
		\end{align*}
		Therefore, by Hölder's inequality, recalling the notation
		\begin{equation*}
			\rho_{2k} := \sum \lambda_j \lt|\opp^{k}\psi_j \rt|^2,
		\end{equation*}
		we get the following bound
		\begin{align}\nonumber
			\|\opp_\ii^{m}(E_\ii)\|_{L^\infty} &\leq \sum_{l=0}^m \binom{m}{l}\lt\||\nabla K| * \lt(\rho_{2l}^{1/2}\rho_{2(m-l)}^{1/2}\rt)\rt\|_{L^\infty}
			\\\label{eq:pm_E_0}
			&\leq C_K \sum_{l=0}^m \binom{m}{l} \|\rho_{2l}\|_{L^{q_1}}^{1/2}\|\rho_{2(m-l)}\|_{L^{q_2}}^{1/2},
		\end{align}
		where $C_K = \|\nabla K\|_{L^\bb}$ and
		\begin{equation}\label{eq:rel_q_q_b}
			\frac{1}{q_1'} + \frac{1}{q_2'} = \frac{2}{\bb}.
		\end{equation}
		From the hypothesis~\eqref{eq:hyp_n} for $n$, we get $n<n_1+1$ and
		\begin{align*}
			(n_1+1-n)\bb \geq \lt(n_1r'+d\rt).
		\end{align*}
		By defining $p'_{n_1,k} := (n_1/k)'p'_{n_1} = (n_1/k)'\lt(r'+d/n_1\rt)$, it implies that
		\begin{equation*}
			\bb \geq p'_{n_1,n-1}.
		\end{equation*}
		Moreover, by the interpolation inequalities \eqref{eq:Lieb_Thirring_k} and the fact that $M_{n_1}$ and $M_0$ are bounded on $[0,T]$, we deduce that $\|\rho_{k}\|_{L^p}$ is bounded uniformly with respect to $\hbar$ for any $t\in[0,T]$ and any $p\in [1,p_{n_1,k}]$. In particular, since
		\begin{equation*}
			\frac{2}{\bb}\leq  \frac{2}{p'_{n_1,n-1}} \leq \frac{2}{p'_{n_1,m}} =	\frac{1}{p'_{n_1,2l}} + \frac{1}{p'_{n_1,2(m-l)}},
		\end{equation*}
		we can find $q_1,q_2\geq 1$ such that the left hand side of \eqref{eq:pm_E_0} is bounded on $[0,T]$ and \eqref{eq:rel_q_q_b} is verified, and there exists $(\eps_1,\eps_2)\in(0,1)^2$ such that
		\begin{align}\nonumber
			\frac{1}{q'_1} &= \frac{\eps_1}{p'_{n_1,2l}}
			\\\nonumber
			\frac{1}{q'_2} &= \frac{\eps_2}{p'_{n_1,2(m-l)}}
			\\\nonumber
			\|\rho_{2l}\|_{L^{q_1}}\|\rho_{2(m-l)}\|_{L^{q_2}} &\leq \|\rho_{2l}\|_{L^{p_{n_1,2l}}}^{\eps_1}\|\rho_{2(m-l)}\|_{L^{p_{n_1,2(m-l)}}}^{\eps_2} \|\rho_{2l}\|_{L^1}^{1-\eps_2}\|\rho_{2(m-l)}\|_{L^1}^{1-\eps_2}
			\\\label{eq:58}
			&\leq C_{d,r,n_1}^2 \|\op\|_{\L^r}^{2\Theta_0} M_{n_1}^{2\Theta_1} M_{2l}^{1-\eps_1} M_{2(m-l)}^{1-\eps_2},
		\end{align}
		where
		\begin{align*}
			\Theta_0 &= \frac{1}{2}\lt(\eps_1 \lt(\frac{r'}{p'_{n_1,2l}}\rt) + \eps_2 \lt(\frac{r'}{p'_{n_1,2(m-l)}}\rt)\rt)
			\\
			\Theta_1 &= \frac{1}{2}\lt(\eps_1 \lt(1-\frac{r'}{p'_{n_1,2l}}\rt) + \eps_2 \lt(1-\frac{r'}{p'_{n_1,2(m-l)}}\rt)\rt).
		\end{align*}
		Since by \eqref{eq:rel_q_q_b}, $\frac{2}{\bb} = \frac{\eps_1}{p'_{n_1,2l}} + \frac{\eps_2}{p'_{n_1,2(m-l)}}$, we deduce that
		\begin{align*}
			\Theta_0 &= \frac{r'}{\bb}
			\\
			\Theta_1 &= \frac{1}{2}\lt(\eps_1+\eps_2\rt)-\frac{r'}{\bb}.
		\end{align*}
		Moreover, by interpolation, for any $k\in[0,n_1]$,
		\begin{equation*}
			M_k \leq M_{n_1}^{k/n_1}M_0^{1-k/n_1} \leq M_0 + M_{n_1}.
		\end{equation*}
		Using this inequality for $k=2l$ and $k=2(m-l)$ in \eqref{eq:58}, inequality \eqref{eq:pm_E_0} becomes
		\begin{equation}\label{eq:pm_E}
			\|\opp_\ii^{m}(E_\ii)\|_{L^\infty} \leq 2^m C_{\op^\init} \lt(1 + M_{n_1}^\theta\rt),
		\end{equation}
		where $\theta = 1-\frac{r'}{\bb}$, $C_{\op^\init} = C_{d,r,n_1}C_K\|\op^\init\|_{\L^r}^{\Theta_0}(1+M_0)$ and we used the propagation of the $\L^r$ and $\L^1$ norm (Proposition~\ref{prop:propag_Lp}).
		We can now come back to \eqref{eq:d_t_lp_m_2}. By combining it with \eqref{eq:util_leibniz}, \eqref{eq:util_holder} and \eqref{eq:pm_E}, we arrive at
		\begin{align*}
			\dt\lt(\lt\| \opp_\ii^n \op \rt\|_{2p}\rt) &= \frac{1}{2p\lt\| \opp_\ii^n \op \rt\|_{2p}^{2p-1}}\dt \Tr(\lt| \opp_\ii^n \op \rt|^{2p})
			\\
			&\leq C_{\op^\init} \lt(1 + M_{n_1}^\theta\rt) \sum_{k=1}^{n}\sum_{m=0}^{n-k}\binom{n-1}{\lfloor (k-1)/2\rfloor}\binom{n-k}{m} 2^m \lt\|\op\rt\|_{2p}^\frac{m+1}{n} \lt\|\opp_\ii^{n}\op\rt\|_{2p}^{1-\frac{m+1}{n}}
			\\
			&\leq 4^n C_{\op^\init} \lt(1 + M_{n_1}^\theta\rt) \lt(\lt\|\op\rt\|_{2p}^\frac{1}{n} \lt\|\opp_\ii^{n}\op\rt\|_{2p}^{1-\frac{1}{n}} + \lt\|\op\rt\|_{2p}\rt).
		\end{align*}
		By Multiplying the inequality by $h^{-d/(2p)'}$ and by conservation of the $\L^{2p}$ norm, we get
		\begin{align*}
			\dt\lt\|\opp_\ii^n \op \rt\|_{\L^{2p}} &\leq 4^n C_{\op^\init} \lt(1 + M_{n_1}^\theta\rt) \lt(\lt\|\op^\init\rt\|_{\L^{2p}}^\frac{1}{n} \lt\|\opp_\ii^{n}\op\rt\|_{\L^{2p}}^{1-\frac{1}{n}} + \lt\|\op^\init\rt\|_{\L^{2p}}\rt).
		\end{align*}
		Defining $u := \frac{\lt\|\opp_\ii^n \op \rt\|_{\L^{2p}}} {\lt\|\op^\init\rt\|_{\L^{2p}}}$ and $c(t) := 4^n C_{\op^\init} \int_0^t\lt(1 + M_{n_1}^\theta\rt)$, it can be written
		\begin{align*}
			\ddt{u} &\leq \lt(1+u^{1-\frac{1}{n}}\rt) \ddt c.
		\end{align*}
		By Gronwall's Lemma, we obtain
		\begin{align*}
			u(t) &\leq (2c(t)+u(0)) + (u(0)^\frac{1}{n} + 2c(t)/n)^n
			\\
			&\leq 2^n (u(0) + c(t)^n),
		\end{align*}
		or equivalently
		\begin{equation}\label{eq:estim_Lp_m}
			\lt\|\opp_\ii^n \op \rt\|_{\L^{2p}} \leq 2^n \lt(\lt\|\opp_\ii^n \op^\init \rt\|_{\L^{2p}} + \tilde{C}_{\op^\init} \lt(t+\int_0^tM_{n_1}^\theta\rt)^n\rt),
		\end{equation}
		where $\tilde{C}_{\op^\init} = \lt\|\op^\init\rt\|_{\L^{2p}} (4^nC_{\op^\init})^n$. It proves inequality \eqref{eq:propag_Lp_m}. Remark that if $r=\infty$, then we can take $C_{\op^\init}$ depending only on $\op^\init$ and not on $p$ since by interpolation between $\L^p$ spaces (Proposition~\ref{eq:Holder_op}), we have
		\begin{equation*}
			\lt\|\op^\init\rt\|_{\L^{2p}} \leq \lt\|\op^\init\rt\|_{\L^{\infty}}^{1/(2p)'}\lt\|\op^\init\rt\|_{\L^{1}}^{1/(2p)}  \leq \lt\|\op^\init\rt\|_{\L^{\infty}} + \lt\|\op^\init\rt\|_{\L^{1}}.
		\end{equation*}
		Therefore we can pass to the limit $p\to\infty$ in \eqref{eq:estim_Lp_m} to get
		\begin{equation*}
			\lt\|\opp_\ii^n \op \rt\|_{\L^\infty} \leq 2^n \lt(\lt\|\opp_\ii^n \op^\init \rt\|_{\L^\infty} + \tilde{C}_{\op^\init} \lt(t+\int_0^tM_{n_1}^\theta\rt)^n\rt),
		\end{equation*}
		with $\tilde{C}_{\op^\init} = 4^{n^2}C_{d,n_1}^n \|\nabla K\|_{L^\bb}^n \|\op^\init\|_{\L^\infty}^{1+n/\bb}(1+M_0)^n$.
	\end{demo}

\section{The quantum coupling estimate}\label{sec:coupling}

	Following the ideas of Loeper in \cite{loeper_uniqueness_2006}, we use the property of displacement convexity of the interpolation between probability measures induced by the optimal transport to deduce the following bound in Wasserstein distance.
	
	\begin{prop}\label{prop:ineq_interpol}
		Let $p\in[1,+\infty]$ and $(\rho_0,\rho_1)\in(L^p\cap\P(\R^d))^2$. Then
		\begin{equation}
			\|\rho_0-\rho_1\|_{W^{-1,\frac{2p}{p+1}}} \leq \max(\|\rho_0\|_{L^p},\|\rho_1\|_{L^p})^\frac{1}{2} W_2(\rho_0,\rho_1),
		\end{equation} 
		where $\dot{W}^{-1,r}$ denotes the dual space of the space
		\begin{equation*}
			\dot{W}^{1,r'} := \lt\{\varphi, \nabla\varphi\in L^{r'}, \varphi\underset{|x|\to\infty}{\longrightarrow} 0\rt\}.
		\end{equation*}
	\end{prop}
	
	\begin{demo}
		Let $q = p'$ be the Hölder conjugate of $p$, $T$ be the optimal transport map for the $W_2$ distance and $\varphi\in \dot{W}^{1,2q}$. Then the interpolant $\rho_\theta = ((1-\theta)x + \theta T(x))_\#\rho_0$ verifies
		\begin{align*}
			\intd \varphi \rho_\theta = \intd \varphi(x_\theta) \rho_0(\!\d x),
		\end{align*}
		where we denote by $x_\theta := (1-\theta)x + \theta T(x)$. By differentiating with respect to $\theta$ and using Cauchy-Schwartz inequality, we get
		\begin{align*}
			\frac{\!\d}{\!\d\theta}\intd \varphi \rho_\theta &= \intd (T(x)-x)\cdot\nabla\varphi(x_\theta) \rho_0(\!\d x)
			\\
			&\leq \lt(\intd |T(x)-x|^2\rho_0(\!\d x)\rt)^\frac{1}{2}\lt(\intd |\nabla\varphi|^2\rho_\theta\rt)^\frac{1}{2}.
		\end{align*}
		The first integral is nothing but the $W_2$ distance between $\rho_0$ and $\rho_1$. Thus, using Hölder's inequality to bound the second integral, we get
		\begin{align*}
			\frac{\!\d}{\!\d\theta}\intd \varphi \rho_\theta &\leq W_2(\rho_0,\rho_1) \|\varphi\|_{\dot{W}^{1,2q}}\|\rho_\theta\|_{L^p}^{1/2}.
		\end{align*}
		By displacement convexity (see for example \cite[Proposition 7.29]{santambrogio_optimal_2015}), the following inequality holds 
		\begin{align*}
			\|\rho_\theta\|_{L^p} \leq \max(\|\rho_0\|_{L^p},\|\rho_1\|_{L^p}).
		\end{align*}
		Noticing that $(2q)' = \frac{2p}{p+1}$, an integration with respect to $\theta$ on $[0,1]$ gives the expected result.
	\end{demo}

	As a consequence of Proposition \ref{prop:ineq_interpol} and the weak Young inequality, we get the following inequality
	\begin{cor}\label{cor:ineq_inertpol}
		Let $p\in(1,+\infty]$, $s=(2p)'$ and $K$ be such that $\nabla^2 K \in L^{s,\infty}$. Then, we have
		\begin{align}
			\|\nabla K*(\rho_0-\rho_1)\|_{L^2} &\leq \|\nabla^2 K\|_{L^{s,\infty}} \max(\|\rho_0\|_{L^p},\|\rho_1\|_{L^p})^\frac{1}{2} W_2(\rho_0,\rho_1).
		\end{align}
		If $p=1$, the same formula holds by replacing $L^{2,\infty}$ by $L^2$ and if $p=\infty$ by replacing $L^{1,\infty}$ by $L^1$.\\
		Moreover, if $p=\infty$, $\|\nabla^2 K\|_{L^1}$ can be replaced by $\|\nabla K\|_{B^1_{1,\infty}}$.
	\end{cor}
	
	\begin{demo}
		Let $r=\frac{2p}{p+1}$. We first write that for $\rho := \rho_0-\rho_1$ and $\varphi\in L^2$,
		\begin{equation}\label{eq:tmp1}
			\lt|\intd (\nabla K*\rho) \varphi\rt|
			\leq \|\rho\|_{\dot{W}^{-1,r}}\|\nabla K* \varphi\|_{\dot{W}^{1,r'}}.
		\end{equation}
		Then, as a consequence of the weak Young inequality (see \cite[Chapter 4, (7)]{lieb_analysis_2001}, we have
		\begin{equation}\label{eq:tmp2}
			\|\nabla K* \varphi\|_{\dot{W}^{1,r'}} = \|\nabla^2 K* \varphi\|_{L^{r'}}
			\leq \|\varphi\|_{L^2}\|\nabla^2 K\|_{L^{s,\infty}},
		\end{equation}
		with $\frac{1}{s} = 1 + \frac{1}{r'} - \frac{1}{2} = 1 - \frac{1}{2p}$. Combining \eqref{eq:tmp1} and \eqref{eq:tmp2}, by duality, we deduce
		\begin{equation*}
			\|\nabla K*\rho\|_{L^2}
			\leq \|\rho\|_{W^{-1,r}}\|\nabla^2 K\|_{L^{s,\infty}}.
		\end{equation*}
		We then use Proposition \ref{prop:ineq_interpol} to conclude. When $p=\infty$ and $r=2$, we use the fact that
		\begin{equation}
			\|\nabla K* \varphi\|_{\dot{H}^1} \leq \|\varphi\|_{L^2}\|\nabla K\|_{B^1_{1,\infty}},
		\end{equation}
		which is proved in Appendix (see \eqref{eq:Besov_2} in Proposition \ref{prop:Besov}).
	\end{demo}

	We can now prove the following key estimate in the modified Wasserstein distance as defined by \eqref{def:Wh}.
		
	\begin{prop}\label{prop:coupling_estimate}
		Let $(s,q)\in(1,2)\times[1,\infty]$ and assume
		\begin{align*}
			\nabla^2 K &\in L^{s,\infty}\cap L^q,
		\end{align*}
		with $L^{s,\infty}$ replaced by $L^2$ if $s=2$. Let $\oph\in\PP$ be a solution of \eqref{eq:Hartree} equation and $f$ be a solution of the \eqref{eq:Vlasov} equation such that the spatial densities verify
		\begin{align*}
			\rho_\hbar := \intd \fht\d\xi &\in L^\infty([0,T],L^{q'}\cap L^{{s'}/2})
			\\
			\rho := \intd f\d \xi &\in L^\infty([0,T],L^\infty),
		\end{align*}
		uniformly with respect to $\hbar$. Then, for all $t \in[0,T]$, we have
		\begin{equation*}
			\Wh(f(t),\oph(t)) \leq \Wh(f^\init,\oph^\init) e^{Ct} + C_0(t)\sqrt{\hbar},
		\end{equation*}
		where
		\begin{align*}
			C_1 &= \|\nabla^2 K\|_{L^{s,\infty} + L^2} \sup_{[0,T]}\lt(\max(\|\rho\|_{L^{s'/2}}, \|\rho_\hbar\|_{L^{s'/2}})^{1/2} \|\rho\|_{L^\infty}^{1/2}\rt)
			\\
			C &= 1 + C_1 + \sup_{[0,T]}\|\rho_\hbar\|_{L^{q'}} \|\nabla^2 K\|_{L^{q}}
			\\
			C_0(t) &= C_1\sqrt{d}(C)^{-1}(e^{Ct}-1).
		\end{align*}
	\end{prop}
	
	\begin{demo}
		Let $p = q'$ and $\tilde{p} = s'/2$. As in \cite[Section 4]{golse_schrodinger_2017}, we define the time dependent coupling $\gam(z) = \gam_\hbar(t,z)$ with $z=(x,\xi)$ as the solution to the Cauchy problem 
		\begin{equation*}
			\partial_t\gam = \{H,\gam\} + \frac{1}{i\hbar}[\Hh,\gam],
		\end{equation*}
		with initial condition $\gam^\init \in \mathcal{C}(f^\init,\oph^\init)$. As proved in \cite[Lemma 4.2]{golse_schrodinger_2017}, $\gam \in \mathcal{C}(f(t),\oph(t))$. We also define
		\begin{equation*}
			\Eh = \Eh(t) := \intdd \Tr\lt(\ch(z)\gam(z)\rt)\d z.
		\end{equation*}
		By differentiating in time, we get
		\begin{equation*}
			\ddt{\Eh} = \intdd \Tr\lt(\lt(\{H,\ch(z)\}+\frac{1}{i\hbar}[\Hh,\ch]\rt)\gam(z)\rt)\d z,
		\end{equation*}
		which, by a direct computation, as detailed in \cite[Section 4.3]{golse_schrodinger_2017}, leads to
		\begin{align}\label{eq:estim_E_0}
			\ddt{\Eh} \leq \Eh &+ \intdd \Tr_y((\xi-\opp)\cdot(E_\hbar(y)-E(x))\gam(z))\d z
			\\\nonumber
			&+ \intdd \Tr_y((E_\hbar(y)-E(x))\cdot(\xi-\opp)\gam(z))\d z.
		\end{align}
		Since $\gam\geq 0$, we use the fact that by Hölder's inequality for Schatten spaces (see e.g. \cite{simon_trace_2005}) and cyclicity of the trace, we have for any operators $(A,B) \in \L(L^2,L^2(\R^d,\R^d))^2$
		\begin{align*}
			\Tr(A^* B\gam)^2 &= \Tr(\gam^{1/2}A^* B\gam^{1/2})^2
			\\
			&\leq \Tr(|\gam^{1/2}A^*|^2)\Tr(|B\gam^{1/2}|^2)
			\\
			&\leq \Tr(A\gam A^*) \Tr(\gam^{1/2}B^*B\gam^{1/2})
			\\
			&\leq \Tr(|A|^2\gam) \Tr(|B|^2\gam).
		\end{align*}
		Thus, using this inequality for $A = (\xi-\opp)$ and $B = E_\hbar(y)-E(x)$ for the first integral in \eqref{eq:estim_E_0} and $A = E_\hbar(y)-E(x)$ and $B = (\xi-\opp)$ for the second integral, we get by Cauchy-Schwartz inequality
		\begin{align}\label{eq:estim_E_1}
			\ddt{\Eh} \leq \Eh &+ 2\lt(\intdd \Tr(|\xi-\opp|^2\gam(z))\d z\rt)^\frac{1}{2} \lt(\intdd \Tr(|E_\hbar(y)-E(x)|^2\gam(z))\d z\rt)^\frac{1}{2}.
		\end{align}
		The first integral is bounded by $\Eh$ and second integral by $2 (I_1 + I_2)$ where
		\begin{align*}
			I_1 &= \intdd \Tr(|E_\hbar(x)-E(x)|^2\gam(z))\d z
			\\
			I_2 &= \intdd \Tr(|E_\hbar(y)-E_\hbar(x)|^2\gam(z))\d z.
		\end{align*}
		Then, since $\gam\in\mathcal{C}(f,\oph)$, by corollary \ref{cor:ineq_inertpol}, we can control $I_1$ in the following way
		\begin{align*}
			I_1 &= \intd |\nabla K*(\rho_\hbar-\rho)(x)|^2\rho(x)\d x
			\\
			&\leq \|\nabla^2 K\|_{L^{s,\infty}}^2 \max(\|\rho\|_{L^{\tilde{p}}},\|\rho_\hbar\|_{L^{\tilde{p}}})  W_2(\rho,\rho_\hbar)^2\ \|\rho\|_{L^\infty}.
		\end{align*}
		Moreover, since $\rho_\hbar = \intd \fht(t,x,\xi)\d\xi$ is nothing but the projection of $f_\hbar$ on the space of positions, we have $W_2(\rho,\rho_\hbar) \leq W_2(f,\fht)$ (see Proposition \ref{prop:projection_wasserstein} for a more detailed proof). Using Theorem~\ref{th:comparaison_W2} and the definition of $\Wh$, we get
		\begin{align}
			I_1 & \leq C_1^2 W_2(f,\fht)^2
			\nonumber\\
			& \leq C_1^2(\Wh(f,\oph)^2 + d\hbar)
			\nonumber\\
			& \leq C_1^2(\Eh + d\hbar).
			\label{eq:estim_I1}
		\end{align}
		In order to control $I_2$, we remark that, from Young's inequality, we get
		\begin{align*}
			\|\nabla E_\hbar\|_{L^\infty} = \|\nabla^2 K* \rho_\hbar\|_{L^\infty} & \leq \|\rho_\hbar\|_{L^p}\|\nabla^2 K\|_{L^q},
		\end{align*}
		which implies that $E_\hbar\in C^{0,1}$ uniformly with respect to $\hbar$, and
		\begin{align*}
			I_2 &\leq C_2^2 \intdd \Tr(|y-x|^2\gh(t,z))\d z \leq C_2^2 \Eh,
		\end{align*}
		where $C_2 = \|\rho_\hbar\|_{L^p} \|\nabla^2 K\|_{L^q}$. By combining this estimate with \eqref{eq:estim_I1}, equation \eqref{eq:estim_E_1} becomes
		\begin{align*}
			\ddt{\Eh} &\leq \Eh + \sqrt{\Eh} \lt(2(C_1^2+C_2^2)\Eh + 2d\hbar C_1^2 \rt)^\frac{1}{2}
			\\
			&\leq \lt(1+\sqrt{2}(C_1+C_2)\rt)\Eh + \sqrt{2d\hbar}C_1 \sqrt{\Eh},
		\end{align*}
		which leads to 
		\begin{equation*}
			\ddt{\sqrt{\Eh}} \leq \lt(1+C_1+C_2\rt)\sqrt{\Eh} + \sqrt{d\hbar}C_1 \hbar. 
		\end{equation*}
		By Grönwall's inequality, it leads to
		\begin{equation*}
			\Wh(f,\oph) \leq \sqrt{\Eh} \leq \sqrt{\Eh}(0) e^{Ct} + C_1 \sqrt{d\hbar}\ \frac{e^{Ct}-1}{C}.
		\end{equation*}
		Minimizing the right hand side as $\gh^\init$ runs through $\mathcal{C}(f^\init,\gh^\init)$ gives the expected result.
	\end{demo}
	
	When $\nabla K\in B^1_{1,\infty}$, which includes the Coulomb potential, previous proposition becomes
	\begin{prop}\label{prop:coupling_estimate_Coulomb}
		Assume 
		\begin{equation*}
			\nabla K\in B^1_{1,\infty}.
		\end{equation*}
		Let $\oph\in\PP$ be a solution of \eqref{eq:Hartree} equation and $f$ be a solution of the \eqref{eq:Vlasov} equation such that the respective spatial densities verify
		\begin{align*}
			\rho_\hbar &\in L^\infty([0,T],L^\infty)
			\\
			\rho &\in L^\infty([0,T],L^\infty),
		\end{align*}
		uniformly with respect to $\hbar$. Then, for all $t \in[0,T]$, we have
		\begin{equation*}
			\Wh(f(t),\oph(t)) \leq \max\lt(\sqrt{d\hbar},\,\Wh(f^\init,\oph^\init)^{e^{t/\sqrt{2}}} e^{\lambda(e^{t/\sqrt{2}}-1)}\rt),
		\end{equation*}
		where
		\begin{align*}
			\lambda &=  C\lt(1 + \|\nabla K\|_{B^1_{1,\infty}}\sup_{[0,T]}(\|\rho\|_{L^\infty} + \|\rho_\hbar\|_{L^\infty})\rt).
		\end{align*}
	\end{prop}
	
	\begin{demo}
		The proof is similar to the proof of Proposition~\ref{prop:coupling_estimate}. With the same notations, we arrive at
		\begin{equation}\label{eq:estim_E_2}
			\ddt{\Eh} \leq \Eh + \sqrt{2\Eh} (I_1 + I_2)^{1/2}.
		\end{equation}
		Then by corollary \ref{cor:ineq_inertpol}, we obtain
		\begin{equation*}
			I_1 \leq \|\nabla K\|_{B^1_{1,\infty}}^2 \max(\|\rho\|_{L^\infty},\|\rho_\hbar\|_{L^\infty})  W_2(\rho,\rho_\hbar)^2\ \|\rho\|_{L^\infty}.
		\end{equation*}
		As in the proof of Proposition~\ref{prop:coupling_estimate}, it leads to
		\begin{equation*}
			I_1 \leq C_1^2(\Eh + d\hbar),
		\end{equation*}
		where $C_1 = \|\nabla K\|_{B^1_{1,\infty}} \max(\|\rho\|_{L^\infty},\|\rho_\hbar\|_{L^\infty})^{1/2} \|\rho\|_{L^\infty}^{1/2}$.
		In order to control $I_2$, we use the fact that since $\nabla K\in B^1_{1,\infty}$, then, as proved in Appendix~\ref{appendix:Besov} (inequality~\eqref{eq:Besov_1} of Proposition~\ref{prop:Besov}), we have
		\begin{equation}\label{eq:estim_I1_2}
			\|E_\hbar\|_{B^1_{\infty,\infty}} = \|\nabla K* \rho_\hbar\|_{B^1_{\infty,\infty}} \leq \|\rho_\hbar\|_{L^\infty}\|\nabla K\|_{B^1_{1,\infty}}.
		\end{equation}
		Then we use a result proved for example in \cite[Chapter 2]{bahouri_fourier_2011} which states that any function in $B^1_{\infty,\infty}$ is log-Lipschitz in the sense that for any $|x-y|<1$, we have
		\begin{align*}
			|E_\hbar(x)-E_\hbar(y)| \leq \|E_\hbar\|_{B^1_{\infty,\infty}} |x-y| \lt(1+\lt|\ln(|x-y|)\rt|\rt).
		\end{align*}
		But for any $r\in(0,1)$, since $B^1_{1,\infty} \subsetArrow L^\infty$, for any $|x-y|\geq r$, we get
		\begin{align*}
			|E_\hbar(x)-E_\hbar(y)| \leq 2\|E_\hbar\|_{L^\infty} \leq C\|E_\hbar\|_{B^1_{\infty,\infty}}\frac{|x-y|}{r}.
		\end{align*}
		Let introduce the kernel of $\gh$, $\gamma(y_1,y_2,z)$ (which still depends on $t$ and $\hbar$) and its diagonal
		\begin{equation*}
			\gamma(y,z) := \gamma(y,y,z).
		\end{equation*}
		Then, we have
		\begin{align*}
			I_2 &= \intdd \intd |E_\hbar(y)-E_\hbar(x)|^2\gamma(y,z) \d y\d z.
			\\
			&\leq C_2^2 \lt(\Eh + \intdd\int_{|x-y|<r} |y-x|^2\ln(|x-y|)^2\gamma(y,z)\d y\d z\rt)
			\\
			&\leq C_2^2 \lt(\Eh + \frac{1}{4}\intdd\int_{|x-y|<r} F(|y-x|^2)\gamma(y,z)\d y\d z\rt),
		\end{align*}
		where $C_2 = \lt(\frac{C}{r} + 1\rt)^{1/2} \|\rho_\hbar\|_{L^\infty} \|\nabla K\|_{B^1_{1,\infty}}$ and $F(x) = x\ln(x)^2$. As noticed in \cite{loeper_uniqueness_2006}, $F$ is concave on $[0,e^{-1}]$. Thus, by taking $r=e^{-1}$, by Jensen's inequality,
		\begin{align*}
			I_2 &\leq C_2^2 \lt(\Eh + \frac{1}{4} F(\Eh)\rt).
		\end{align*}
		By combining this estimate with \eqref{eq:estim_I1_2}, equation \eqref{eq:estim_E_2} becomes
		\begin{align*}
			\ddt{\Eh} &\leq \Eh + \sqrt{2\Eh} ((C_1^2+C_2^2)\Eh + C_1^2 d\hbar + F(\Eh)/4)^{1/2}
			\\
			&\leq (1+\sqrt{2}(C_1+C_2))\Eh + C_1\sqrt{2d\hbar\Eh} + \Eh \ln(\Eh)/\sqrt{2}
			\\
			&\leq \lambda\Eh + C_1\sqrt{d\hbar/2} + \Eh \ln(\Eh)/\sqrt{2},
		\end{align*}
		where $\lambda = 1+\sqrt{2}(2C_1+C_2)$ and we used the inequalities $\sqrt{a+b} \leq \sqrt{a} + \sqrt{b}$ and $\sqrt{2ab}\leq a+b$.
		Then, for any $t$ such that $\lambda\Eh\geq C_1\sqrt{d\hbar/2}$, we get
		\begin{align*}
			\ddt{\ln(\Eh)} &\leq 2\lambda + \ln(\Eh)/\sqrt{2}.
		\end{align*}
		By Grönwall's inequality, it leads to
		\begin{equation*}
			\Wh(f,\oph) \leq \Eh \leq \max\lt(\frac{C_1(t)\sqrt{d\hbar}}{\lambda(t)\sqrt{2}},\,\Eh(0)^{e^{t/\sqrt{2}}} e^{\sqrt{2}\tilde{\lambda}(t)(e^{t/\sqrt{2}}-1)}\rt),
		\end{equation*}
		where $\tilde{\lambda}(t) = \sup_{[0,T]}\lambda$, and which gives the expected result since $C_1(t) \leq \lambda(t)$.
	\end{demo}
	
	Combining the propagation of moments of Theorem~\ref{th:propag_moments} with the Proposition~\ref{prop:coupling_estimate} which gives the semiclassical convergence as soon as $\rho$ is sufficiently integrable, we can now prove Theorem~\ref{th:CV}. Theorem~\ref{th:CV_VP} is proved in the same way using Proposition~\ref{prop:coupling_estimate_Coulomb} and Proposition~\ref{prop:propag_Lp_m}.
	
	\begin{demo}[ of Theorem~\ref{th:CV}]
		Since $\nabla K \in L^\infty + L^{\bb,\infty}$, from Theorem~\ref{th:propag_moments}, we obtain the existence of $T\in(0,+\infty]$ and $\Phi\in C^0([0,T))$ such that for any $t\in[0,T)$
		\begin{equation*}
			M_{n_1} < \Phi(t).
		\end{equation*}
		Moreover, from Proposition~\ref{prop:propag_Lp} we know that
		\begin{equation*}
			\|\oph\|_{\L^r} = \|\oph^\init\|_{\L^r} \leq C.
		\end{equation*}
		By inequality~\eqref{eq:Lieb_Thirring_+}, we deduce that 
		\begin{equation*}
			\|\rho_\hbar\|_{L^{p_{n_1}}} \leq \Phi(t)^{1-\theta}.
		\end{equation*}
		Moreover, by Proposition~\ref{prop:propag_Lp}, we also deduce the propagation of the mass
		\begin{equation*}
			\intd \rho_\hbar = \Tr(\oph) = \|\oph\|_{\L^1} = \|\oph^\init\|_{\L^1} = 1.
		\end{equation*}
		Remarking that
		\begin{align*}
			q \geq p'_{n_1} &\ssi q \geq r' + \frac{d}{n_1} \ssi n_1 \geq \frac{d}{q-r'},
		\end{align*}
		we get that $p := q' \in [1,p'_{n_1}]$. Moreover, since $q'\geq2$, it also implies that $q'/2\in [1,p_{n_1}']$. By Hölder's inequality, it implies that for a given $\eps < 1-\theta$,
		\begin{align*}
			\|\rho_\hbar\|_{L^{q'}} &\leq \Phi(t)^{\eps}
			\\
			\|\rho_\hbar\|_{L^{q'/2}} &\leq \Phi(t)^{2\eps},
		\end{align*}
		and  we can use Proposition~\ref{prop:coupling_estimate}  to get the result.
	\end{demo}

\section{Superpositions of coherent states}\label{sec:coh_stat}

	We recall in this section some results about the approximation of measures on the phase space by a superposition of coherent states and state some applications in our case. See also Thirring \cite{thirring_quantum_1983}, Lions and Paul \cite{lions_sur_1993}, Golse et al \cite{golse_mean_2016}. Let $\varphi\in L^1$ be a smooth function such that $\|\varphi\|_{L^2} = 1$. Then the coherent states are defined by
	\begin{equation*}
		\varphi_{x,\xi}(y) = \frac{1}{h^{d/4}}\varphi\lt(\frac{y-x}{\sqrt{h}}\rt) e^{2i\pi y\cdot\xi/h},
	\end{equation*}
	and we will denote the associated density operator by
	\begin{equation*}
		\op_{x,\xi} := \ket{\varphi_{x,\xi}}\bra{\varphi_{x,\xi}}.
	\end{equation*}
	We can then associate to a measure $\mu\in\P(\R^{2d})$ of the phase space the following operator
	\begin{equation*}
		\opmu := \OP_\varphi(\mu) := \iintd \op_{x,\xi}\mu(\!\d x\d\xi).
	\end{equation*}
	It corresponds to the density operator defined in \cite[Exemple III.7]{{lions_sur_1993}}. Up to a constant depending on $\hbar$, this is also what is called a Töplitz operator in \cite{golse_mean_2016}. The constant comes from the fact that we consider operators associated to measures with finite mass on the semiclassical limit, while Töplitz operators describe operators acting on these measures.
	
	As expected, the mass is the trace of the operator
	\begin{equation*}
		\iintd \mu = \Tr(\opmu).
	\end{equation*}
	Moreover, we remark that $\op_{x,\xi} = \OP_\varphi(\delta_{x,\xi})$ and as proved in \cite{lions_sur_1993}, by defining the Wigner transform $\delta_{x,\xi}^\varphi := \wh(\op_{x,\xi})$, the following holds
	\begin{align}\nonumber
		\delta_{x,\xi}^\varphi &\underset{h\to 0}{\rightharpoonup} \delta_{x,\xi}
		\\\label{eq:approx_dirac}
		\wh(\opmu) = \delta_{0,0}^\varphi * \mu &\underset{h\to 0}{\rightharpoonup} \mu,
	\end{align}
	where the convergence holds in the sense of the duality with $C_0(\R^{2n})$. An other result proved in \cite{golse_schrodinger_2017} is the comparison between the Wasserstein pseudo-distance defined in \eqref{def:Wh} with the classical Wasserstein pseudo-distance, which completes Theorem~\ref{th:comparaison_W2}
	\begin{prop}[Golse and Paul \cite{golse_schrodinger_2017}]
		Let $(\mu,\nu)\in\P(\R^{2d})^2$ be two probability measures such that $W_2(\nu,\mu) < \infty$ and $\opmu := \OP_\varphi(\mu)$ where $\varphi$ is a Gaussian with $\|\varphi\|_{L^2} = 1$. Then
		\begin{equation*}
			|\Wh(\nu,\opmu) - W_2(\nu,\mu)| \leq \sqrt{2d\hbar}.
		\end{equation*}
	\end{prop}
	
	Finally, the following proposition justifies our definition \eqref{def:q_schatten} for the quantum Lebesgue norm. 
	\begin{prop}
		Let $\mu \in \P(\R^{2d})$ and $\opmu := \OP_\varphi(\mu)$. Then for any $r\geq 1$, it holds
		\begin{align}\label{eq:comparaison_L_r}
			\|\mu\|_{L^r_{x,\xi}} &\leq \|\opmu\|_{\L^r}
			\\\label{eq:comparaison_L_infty}
			\|\mu\|_{L^\infty_{x,\xi}} &= \|\opmu\|_{\L^\infty}.
		\end{align}
		Moreover, in the particular case $\delta_{0,0}^\varphi \geq 0$, we have for any $r\geq 2$
		\begin{equation}\label{eq:comparaison_L_r_2}
			\|\wh(\opmu)\|_{L^r_{x,\xi}} \leq \|\mu\|_{L^r_{x,\xi}} = \|\opmu\|_{\L^r},
		\end{equation}
		with equality in the first inequality if $r=2$, as well as the following convergences
		\begin{align}\label{eq:limit_L_r}
			\|\wh(\opmu)\|_{L^r_{x,\xi}} &\underset{h\to 0}{\rightarrow} \|\opmu\|_{\L^r}
			\\\label{eq:limit_L_r_2}
			\wh(\opmu) &\underset{h\to 0}{\rightarrow} \mu \text{ in } L^r.
		\end{align}
	\end{prop}
	
	\begin{remark}
		The assumption $\delta_{0,0}^\varphi \geq 0$ is verified for example when $\varphi(x) = e^{-\pi|x|^2/2}$, since we can then compute explicitly
		\begin{equation*}
			\delta_{x_0,\xi_0}^\varphi = \frac{1}{h^d} e^{-\frac{\pi}{h}(|x-x_0|^2+|y-y_0|^2)}.
		\end{equation*}
	\end{remark}

	\begin{demo}
		As proved in \cite{thirring_quantum_1983} or \cite[Exemple III.7]{lions_sur_1993}, for any convex mapping $F\geq 0$ such that $F(0) = 0$, it holds
		\begin{equation*}
			\iintd F(\mu) \frac{\!\d x\d\xi}{h^d} \leq \Tr\lt(F\lt(\frac{\opmu}{h^d}\rt)\rt).
		\end{equation*}
		By taking $F(x) = |x|^r$ for $r\geq 1$, it implies in particular
		\begin{equation*}
			\|\mu\|_{L^r_{x,\xi}}^r \leq h^{-d(r-1)} \|\opmu\|_r^r = \|\opmu\|_{\L^r}^r,
		\end{equation*}
		which proves \eqref{eq:comparaison_L_r}. As noticed in \cite[Appendix B]{golse_mean_2016}, this inequality also holds in the other direction when $r=\infty$, which leads to \eqref{eq:comparaison_L_infty}.
		Then, as noticed in \cite{lions_sur_1993}, we deduce from \eqref{eq:approx_dirac} that if $\delta_{0,0}^\varphi \geq 0$, we have
		\begin{equation*}
			\iintd F(\wh(\opmu)) \leq \iintd F(\mu).
		\end{equation*}
		Taking again $F(x)=|x|^r$ leads to the first part of \eqref{eq:comparaison_L_r_2}
		\begin{equation*}
			\|\wh(\opmu)\|_{L^r_{x,\xi}} \leq \|\mu\|_{L^r_{x,\xi}} \leq \|\opmu\|_{\L^r}.
		\end{equation*}
		However, for $r=2$, the following equality holds for any operator $\opmu$
		\begin{equation*}
			\|\wh(\opmu)\|_{L^2_{x,\xi}} = \|\opmu\|_{\L^2}.
		\end{equation*}
		Thus, the above inequalities are equalities when $r=2$
		\begin{equation*}
			\|\wh(\opmu)\|_{L^2_{x,\xi}} = \|\mu\|_{L^2_{x,\xi}} = \|\opmu\|_{\L^2}.
		\end{equation*}
		By complex interpolation, we deduce from the above equation and formula~\eqref{eq:comparaison_L_infty} that for any $r\geq 2$, $\OP_\varphi \in \B(L^r_{x,\xi}, \L^r)$ and
		\begin{equation*}
			\|\opmu\|_{\L^r} \leq \|\mu\|_{L^r_{x,\xi}},
		\end{equation*}
		which proves the equality in formula \eqref{eq:comparaison_L_r_2}. Finally, from \eqref{eq:approx_dirac} and \eqref{eq:comparaison_L_r_2}, we deduce that $\wh(\opmu) \rightharpoonup \mu$ in $L^r$ and
		\begin{equation*}
			\|\mu\|_{L^r_{x,\xi}} \leq \liminf_{h\to 0}\|\wh(\opmu)\|_{L^r_{x,\xi}},
		\end{equation*}
		which combined with \eqref{eq:comparaison_L_r_2} leads to \eqref{eq:limit_L_r} and then \eqref{eq:limit_L_r_2}.
	\end{demo}
	
	Combining all these results, we can for example write a simplified version of Theorem~\ref{th:CV} for superposition of coherent states.
	
	\begin{thm}
		Assume $K$ verifies \eqref{hyp:grad_1} and \eqref{hyp:grad_2} and let $f$ be a solution of the \eqref{eq:Vlasov} equation and $\oph$ be a solution of \eqref{eq:Hartree} equation with respective initial conditions 
		\begin{align*}
			f^\init &\in \P\cap L^\infty_{x,\xi} \text{ verifying \eqref{hyp:VP_moment} and \eqref{hyp:VP_intg}}
			\\
			\oph^\init = \OP_\varphi(g^\init) \text{ with } g^\init &\in \P\cap L^\infty_{x,\xi},
		\end{align*}
		 where $\varphi$ is a normalized Gaussian. Assume also that the initial quantum velocity moment
		\begin{equation}
			M_{n_1}^\init < C \text{ for a given } n_1 \geq \frac{d}{q-1}.
		\end{equation}
		Then there exists $T>0$ such that for any $t\in(0,T)$,
		\begin{equation*}
			W_2(f(t),\fht(t)) \leq C_T \lt(W_2(f^\init,g^\init) + \sqrt{\hbar}\rt),
		\end{equation*}
		where $\fht$ is the Husimi transform of $\oph$.
	\end{thm}
	
	The advantage is that in the above results, the semiclassical estimate is stated only in terms of the classical Wasserstein distance which is a true distance, and it also allows to take $f^\init = g^\init$. We can do the same for Theorem~\ref{th:CV_VP}. We state it here for the Coulomb potential in dimension $d=3$.
	\begin{thm}
		Assume $K = \frac{1}{|x|}$ and let $f$ be a solution of the \eqref{eq:Vlasov} equation and $\oph$ be a solution of \eqref{eq:Hartree} equation with respective initial conditions 
		\begin{align*}
			f^\init &\in \P\cap L^\infty_{x,\xi} \text{ verifying \eqref{hyp:VP_moment} and \eqref{hyp:VP_intg}}
			\\
			\oph^\init = \OP_\varphi(g^\init) \text{ with } g^\init &\in \P\cap L^\infty_{x,\xi},
		\end{align*}
		where $\varphi$ is a normalized Gaussian. Moreover, assume that
		\begin{equation*}
			\forall \ii\in\Int{1,3},\, \opp_\ii^4\oph^\init\in \L^\infty,
		\end{equation*}
		where $\opp_\ii := -i\hbar\partial_\ii$. Assume also that the initial quantum velocity moment
		\begin{equation}
			M_{16}^\init < C.
		\end{equation}
		Then there exists $T>0$ such that
		\begin{align*}
			\rho_\hbar &\in L^\infty([0,T],L^\infty),
		\end{align*}
		uniformly in $\hbar$, and there exists a constant $C_T$ depending only on the intial conditions and independent of $\hbar$ such that
		\begin{equation*}
			W_2(f(t),\fht(t)) \leq C_T \lt(W_2(f^\init,g^\init) + \sqrt{\hbar}\rt).
		\end{equation*}
	\end{thm}

\appendix

\section{Besov Spaces}\label{appendix:Besov}

	We recall that a possible definition of Besov spaces (see e.g.~\cite{bahouri_fourier_2011}) can be done by defining the following norm
	\begin{equation}\label{def:Besov}
		\|u\|_{B^{s}_{p,r}} = \lt\|\lt(2^{sj}\|\Delta_ju\|_{L^p}\rt)_{j\in\mathbb{Z}}\rt\|_{\ell^r},
	\end{equation}
	where $\Delta_j$ is defined by
	\begin{align*}
		\Delta_j u &= 0 &&\text{when } j\leq -2
		\\
		\Delta_{-1} u &= \hat{\chi}*u
		\\
		\Delta_j u &= \F_y(\varphi(2^{-j}y))*u &&\text{when } j\geq 0,
	\end{align*}
	with
	\begin{align}\nonumber
		\chi &\in C^\infty_c(B(0,4/3),[0,1])
		\\\nonumber
		\varphi &\in C^\infty_c\lt(B\lt(0,8/3\rt)\backslash B\lt(0,3/4\rt),[0,1]\rt)
		\\\label{eq:partition}
		\chi &+ \sum_{j\geq 0} \varphi(2^{-j}\cdot) = 1.
	\end{align}
	We also define the space of log-Lipschitz functions by defining the norm
	\begin{equation*}
		\|u\|_{LL} = \sup_{|x-y|\in(0,1)}\lt(\frac{|u(x)-u(y)|}{|x-y|\lt(1+\lt|\ln(|x-y|)\rt|\rt)}\rt),
	\end{equation*}
	for measurable functions $u$ vanishing at infinity. We have the following properties of Besov spaces
	\begin{prop}\label{prop:Besov}
		\begin{equation}\label{eq:Besov_LL}
			B^{1}_{\infty,\infty} \subsetArrow LL.
		\end{equation}
		If $K$ is the Coulomb potential such that $\Delta K = \delta_0$, then we get
		\begin{equation}\label{eq:Besov_Coulomb}
			\lt|\nabla K\rt| = \frac{C}{|x|^{d-1}} \in B^1_{1,\infty}.
		\end{equation}
		If $v\in L^\infty$ and $u\in B^{1}_{1,\infty}$, then
		\begin{align}\label{eq:Besov_1}
			\|u*v\|_{B^1_{\infty,\infty}} &\leq \lt\|u\rt\|_{B^1_{1,\infty}} \|v\|_{L^\infty}
			\\\label{eq:Besov_2}
			\|u*v\|_{\dot{H}^1} &\leq C \|u\|_{B^1_{1,\infty}} \|v\|_{L^2}.
		\end{align}
	\end{prop}
	
	\begin{demo}
		The proof of \eqref{eq:Besov_LL} and \eqref{eq:Besov_Coulomb} can be found for example in \cite[Chapter 2]{bahouri_fourier_2011}.
		To prove \eqref{eq:Besov_1}, we remark that since $\Delta_j$ is a convolution by a smooth and rapidly decaying function, $\Delta_j(u*v) = \Delta_j(u)*v$. By Hölder's inequality, we deduce the following inequality
		\begin{align*}
			\|u*v\|_{B^1_{\infty,\infty}} &= \lt\|\lt(2^{j}\|\Delta_ju*v\|_{L^\infty}\rt)_{j\in\mathbb{Z}}\rt\|_{\ell^\infty}
			\\
			&\leq \|v\|_{L^\infty} \lt\|\lt(2^{j}\|\Delta_ju\|_{L^1}\rt)_{j\in\mathbb{Z}}\rt\|_{\ell^\infty} = \|u\|_{B^1_{1,\infty}} \|v\|_{L^\infty}.
		\end{align*}
		To prove \eqref{eq:Besov_2}, we use the Fourier definition of $\dot{H}^1$ and the fact the Fourier transform is an isometry on $L^2$ to obtain
		\begin{align*}
			\|u*v\|_{\dot{H}^1} &\leq C\||y|\hat{u}(y)\hat{v}(y)\|_{L^2_y} \leq C\||y|\hat{u}(y)\|_{L^\infty_y} \|v\|_{L^2}.
		\end{align*}
		Then by using the fact that $\varphi(2^{-j} y)>0 \ssi |y|\in 2^j[3/4,8/3]$, we obtain the existence of $j_y\geq -2$ such that $\varphi(2^{-j}y) = 0$ for any $j\notin\{j_y-1,j_y,j_y+1\}$ (If $j_y=-2$, then it means that $\chi(y)>0$). Then, by \eqref{eq:partition}, we get
		\begin{align*}
			\|y\hat{u}(y)\|_{L^\infty} &= \lt\|\lt(\chi(y) + \sum_{j\geq 0} \varphi(2^{-j}y)\rt)|y|\hat{u}(y)\rt\|_{L^\infty_y}
			\\
			&\leq C \lt\|\sum_{k=-1}^12^{j_y+k}\F(\Delta_{j_y+k}u)(y)\rt\|_{L^\infty}
			\\
			&\leq C \sup_{j\in\mathbb{Z}}\lt(2^{j} \|\F(\Delta_ju)\|_{L^\infty}\rt)
			\leq C \lt\|\lt(2^{j} \|\Delta_ju\|_{L^1}\rt)_{j\in\mathbb{Z}}\rt\|_{\ell^\infty}.
		\end{align*}
		Therefore, by the definition~\eqref{def:Besov}, we obtain \eqref{eq:Besov_2}.
	\end{demo}

\section{Wasserstein distances}\label{appendix:Wasserstein}

	We recall the definition of the classical Wasserstein-(Monge-Kantorovich) distances between two probability measures $(\mu_0,\mu_1)\in \P(X)^2$ on a given separable Banach space $X$. We first define the notion of coupling by saying that $\gamma \in \P(X^2)$ is a coupling of $\mu_0$ and $\mu_1$ when 
	\begin{equation*}
		(\pi_1)_\#\gamma = \mu_0 \text{ and } (\pi_2)_\#\gamma = \mu_1,
	\end{equation*}
	where $\pi_1$ and $\pi_2$ are respectively the projection on the first and second variable and $\pi_\#\gamma$ denotes the pushforward of the measure $\gamma$ by the map $\pi$. In other words
	\begin{equation*}
		\forall\varphi\in C_0(X), \int_{X^2} \varphi(x) \gamma(\!\d x\d y) = \int_{X} \varphi(x) \mu_0(\!\d x).
	\end{equation*}
	We denote by $\Pi(\mu_0,\mu_1)$ the set of couplings of $\mu_0$ and $\mu_1$. Then we define the Wasserstein-(Monge-Kantorovich) distance in the following way
	\begin{equation}\label{def:Wasserstein}
		W_p(\mu_0,\mu_1) := \lt(\inf_{\gamma\in\Pi(\mu_0,\mu_1)}\int_{X^2} \|x-y\|_X^p\gamma(\!\d x\d y)\rt)^\frac{1}{p}.
	\end{equation}
	The existence of a minimizer is well known and we refer for example to the books \cite{villani_topics_2003} or \cite{santambrogio_optimal_2015} for more properties of these distances.
	
	The following proposition may be classical but we prove it for the sake of completeness
	\begin{prop}\label{prop:projection_wasserstein}
		Let $(f_0,f_1)\in\P(\R^{2d})^2$ and for $\ii\in\{0,1\}$, let $\rho_\ii = (\pi_1)_\# f_\ii$. Then
		\begin{equation*}
			W_2(\rho_0,\rho_1) \leq W_2(f_0,f_1).
		\end{equation*}
	\end{prop}
	
	\begin{demo}
		Let $\gamma\in\P(\R^{2d}\times \R^{2d})$ be the optimal transport plan from $f_0$ to $f_1$ and define $\gamma_\rho = (\pi_{1,3})_\#\gamma$ by
		\begin{equation*}
			\forall\varphi\in C_0(\R^{2d}), \intdd \varphi(x,y)\gamma_\rho(\!\d x\d y) := \int_{\R^{4d}} \varphi(x,y)\gamma(\!\d x\d\xi\d y\d\eta).
		\end{equation*}
		Then for any $\varphi\in C_0$, since the first marginal of $\gamma$ is $f_0$,
		\begin{align*}
			\intdd \varphi(x)\gamma_\rho(\!\d x\d y) &= \int_{\R^{4d}} \varphi(x)\gamma(\!\d x\d\xi\d y\d\eta)
			\\
			&= \intdd \varphi(x)f_0(\!\d x\d\xi)
			\\
			&= \intd \varphi(x)\rho_0(\!\d x).
		\end{align*}
		Hence, the first marginal of $\gamma_\rho$ is $\rho_0$. In the same way, the second marginal of $\gamma_\rho$ is $\rho_1$, and we deduce that $\gamma_\rho\in \Pi(\rho_0,\rho_1)$. Next, let $(\varphi_n)_{n\in\N} \in (C_0(\R^{2d})\cap L^1(\gamma_\rho))^\N$ be an increasing sequence of nonnegative functions converging pointwise to $(x,y)\mapsto|x-y|^2$. By definition of $\gamma_\rho$, for any $n\in\N$, $\varphi\in L^1(\gamma)$. Therefore, by the monotone convergence theorem,
		\begin{align*}
			\intdd |x-y|^2 \gamma_\rho(\!\d x\d y) &= \lim\limits_{n\to\infty} \intdd\varphi_n(x,y)\gamma_\rho(\!\d x\d y)
			\\
			&= \lim\limits_{n\to\infty} \int_{\R^{4d}} \varphi_n(x,y)\gamma(\!\d x\d\xi\d y\d\eta)
			\\
			&= \int_{\R^{4d}} |x-y|^2\gamma(\!\d x\d\xi\d y\d\eta)
			\\
			&\leq \int_{\R^{4d}} (|x-y|^2+|\xi-\eta|^2)\gamma(\!\d x\d\xi\d y\d\eta) = W_2(f_0,f_1)^2.
		\end{align*}
		By definition \eqref{def:Wasserstein}, we deduce
		\begin{equation*}
			W_2(\rho_0,\rho_1)^2 \leq \intdd |x-y|^2 \gamma_\rho(\!\d x\d y) \leq W_2(f_0,f_1)^2,
		\end{equation*}
		which proves the result.
	\end{demo}


\renewcommand{\bibname}{\centerline{Bibliography}}
\bibliographystyle{abbrv} 
\bibliography{Vlasov}

\begin{thebibliography}{10}

\bibitem{ambrosio_semiclassical_2011}
L.~Ambrosio, A.~Figalli, G.~Friesecke, J.~Giannoulis, and T.~Paul.
\newblock Semiclassical limit of quantum dynamics with rough potentials and
  well-posedness of transport equations with measure initial data.
\newblock {\em Communications on Pure and Applied Mathematics},
  64(9):1199--1242, 2011.

\bibitem{ambrosio_passage_2010}
L.~Ambrosio, G.~Friesecke, and J.~Giannoulis.
\newblock Passage from {Quantum} to {Classical} {Molecular} {Dynamics} in the
  {Presence} of {Coulomb} {Interactions}.
\newblock {\em Communications in Partial Differential Equations},
  35(8):1490--1515, July 2010.

\bibitem{amour_classical_2013}
L.~Amour, M.~Khodja, and J.~Nourrigat.
\newblock The classical limit of the {Heisenberg} and time-dependent
  {Hartree}–{Fock} equations: the {Wick} symbol of the solution.
\newblock {\em Mathematical Research Letters}, 20(1):119--139, Jan. 2013.

\bibitem{amour_semiclassical_2013}
L.~Amour, M.~Khodja, and J.~Nourrigat.
\newblock The {Semiclassical} {Limit} of the {Time} {Dependent} {Hartree}–
  {Fock} {Equation}: the {Weyl} {Symbol} of the {Solution}.
\newblock {\em Analysis \& PDE}, 6(7):1649--1674, 2013.

\bibitem{araki_inequality_1990}
H.~Araki.
\newblock On an inequality of {Lieb} and {Thirring}.
\newblock {\em Letters in Mathematical Physics}, 19(2):167--170, 1990.

\bibitem{athanassoulis_strong_2011}
A.~Athanassoulis, T.~Paul, F.~Pezzotti, and M.~Pulvirenti.
\newblock Strong {Semiclassical} {Approximation} of {Wigner} {Functions} for
  the {Hartree} {Dynamics}.
\newblock {\em Rendiconti Lincei - Matematica e Applicazioni}, 22(4):525--552,
  2011.
\newblock arXiv: 1009.0470.

\bibitem{bach_kinetic_2016}
V.~Bach, S.~Breteaux, S.~Petrat, P.~Pickl, and T.~Tzaneteas.
\newblock Kinetic {Energy} {Estimates} for the {Accuracy} of the
  {Time}-{Dependent} {Hartree}-{Fock} {Approximation} with {Coulomb}
  {Interaction}.
\newblock {\em Journal de Mathématiques Pures et Appliquées}, 105(1):1--30,
  2016.

\bibitem{bahouri_fourier_2011}
H.~Bahouri, J.-Y. Chemin, and R.~Danchin.
\newblock {\em Fourier {Analysis} and {Nonlinear} {Partial} {Differential}
  {Equations}}, volume 343 of {\em Grundlehren der mathematischen
  {Wissenschaften}}.
\newblock Springer, Berlin, Heidelberg, Jan. 2011.

\bibitem{bardos_derivation_2002}
C.~Bardos, L.~Erdös, F.~Golse, N.~J. Mauser, and H.-T. Yau.
\newblock Derivation of the {Schrödinger}–{Poisson} {Equation} from the
  {Quantum} {N}-body {Problem}.
\newblock {\em Comptes Rendus Mathematique}, 334(6):515--520, Jan. 2002.

\bibitem{bardos_weak_2000}
C.~Bardos, F.~Golse, and N.~J. Mauser.
\newblock Weak {Coupling} {Limit} of the {N}-particle {Schrödinger}
  {Equation}.
\newblock {\em Methods and Applications of Analysis}, 7(2):275--294, 2000.

\bibitem{benedikter_mean-field_2016}
N.~Benedikter, V.~Jaksic, M.~Porta, C.~Saffirio, and B.~Schlein.
\newblock Mean-{Field} {Evolution} of {Fermionic} {Mixed} {States}.
\newblock {\em Communications on Pure and Applied Mathematics},
  69(12):2250--2303, 2016.

\bibitem{benedikter_hartree_2016}
N.~Benedikter, M.~Porta, C.~Saffirio, and B.~Schlein.
\newblock From the {Hartree} {Dynamics} to the {Vlasov} {Equation}.
\newblock {\em Archive for Rational Mechanics and Analysis}, 221(1):273--334,
  2016.

\bibitem{benedikter_mean-field_2014}
N.~Benedikter, M.~Porta, and B.~Schlein.
\newblock Mean-field {Evolution} of {Fermionic} {Systems}.
\newblock {\em Communications in Mathematical Physics}, 331(3):1087--1131, Nov.
  2014.

\bibitem{brezzi_three-dimensional_1991}
F.~Brezzi and P.~A. Markowich.
\newblock The {Three}-{Dimensional} {Wigner}-{Poisson} {Problem}: {Existence},
  {Uniqueness} and {Approximation}.
\newblock {\em Mathematical Methods in the Applied Sciences}, 14(1):35--61,
  Jan. 1991.

\bibitem{castella_l2_1997}
F.~Castella.
\newblock L2 {Solutions} to the {Schrödinger}–{Poisson} {System}:
  {Existence}, {Uniqueness}, {Time} {Behaviour}, and {Smoothing} {Effects}.
\newblock {\em Mathematical Models and Methods in Applied Sciences},
  7(08):1051--1083, 1997.

\bibitem{egorov_moments_1995}
Y.~V. Egorov and V.~A. Kondratiev.
\newblock On {Moments} of {Negative} {Eigenvalues} of an {Elliptic} {Operator}.
\newblock In {\em Partial {Differential} {Operators} and {Mathematical}
  {Physics}}, Operator {Theory} {Advances} and {Applications}, pages 119--126.
  Birkhäuser Basel, 1995.

\bibitem{erdos_derivation_2001}
L.~Erdös and H.-T. Yau.
\newblock Derivation of the {Nonlinear} {Schrödinger} {Equation} from a {Many}
  {Body} {Coulomb} {System}.
\newblock {\em Advances in Theoretical and Mathematical Physics},
  5(6):1169--1205, 2001.

\bibitem{fournais_semi-classical_2018}
S.~Fournais, M.~Lewin, and J.~P. Solovej.
\newblock The semi-classical limit of large fermionic systems.
\newblock {\em Calculus of Variations and Partial Differential Equations},
  57(4):105, 2018.

\bibitem{frohlich_microscopic_2011}
J.~Fröhlich and A.~Knowles.
\newblock A {Microscopic} {Derivation} of the {Time}-{Dependent}
  {Hartree}-{Fock} {Equation} with {Coulomb} {Two}-{Body} {Interaction}.
\newblock {\em Journal of Statistical Physics}, 145(1):23, Oct. 2011.
\newblock arXiv: 0810.4282.

\bibitem{frohlich_mean-field_2009}
J.~Fröhlich, A.~Knowles, and S.~Schwarz.
\newblock On the {Mean}-{Field} {Limit} of {Bosons} with {Coulomb} {Two}-{Body}
  {Interaction}.
\newblock {\em Communications in Mathematical Physics}, 288(3):1023--1059, June
  2009.
\newblock arXiv: 0805.4299.

\bibitem{gasser_semiclassical_1998}
I.~Gasser, R.~Illner, P.~A. Markowich, and C.~Schmeiser.
\newblock Semiclassical, asymptotics and dispersive effects for
  {Hartree}-{Fock} systems.
\newblock {\em ESAIM: Mathematical Modelling and Numerical Analysis},
  32(6):699--713, 1998.

\bibitem{ginibre_class_1980}
J.~Ginibre and G.~Velo.
\newblock On a {Class} of non {Linear} {Schrödinger} {Equations} with non
  {Local} {Interaction}.
\newblock {\em Mathematische Zeitschrift}, 170(2):109--136, 1980.

\bibitem{ginibre_global_1985}
J.~Ginibre and G.~Velo.
\newblock The {Global} {Cauchy} {Problem} for the non {Linear} {Schrödinger}
  {Equation} {Revisited}.
\newblock {\em Annales de l'Institut Henri Poincare (C) Non Linear Analysis},
  2(4):309--327, July 1985.

\bibitem{golse_mean_2016}
F.~Golse, C.~Mouhot, and T.~Paul.
\newblock On the {Mean} {Field} and {Classical} {Limits} of {Quantum}
  {Mechanics}.
\newblock {\em Communications in Mathematical Physics}, 343(1):165--205, 2016.
\newblock arXiv: 1502.06143.

\bibitem{golse_empirical_2017}
F.~Golse and T.~Paul.
\newblock Empirical {Measures} and {Quantum} {Mechanics}: {Application} to the
  {Mean}-{Field} {Limit}.
\newblock {\em ArXiv e-prints}, 1711:arXiv:1711.08350, Nov. 2017.

\bibitem{golse_schrodinger_2017}
F.~Golse and T.~Paul.
\newblock The {Schrödinger} {Equation} in the {Mean}-{Field} and
  {Semiclassical} {Regime}.
\newblock {\em Archive for Rational Mechanics and Analysis}, 223(1):57--94,
  2017.

\bibitem{golse_wave_2018}
F.~Golse and T.~Paul.
\newblock Wave {Packets} and the {Quadratic} {Monge}-{Kantorovich} {Distance}
  in {Quantum} {Mechanics}.
\newblock {\em Comptes Rendus Mathematique}, 356(2):177--197, Feb. 2018.

\bibitem{golse_derivation_2018}
F.~Golse, T.~Paul, and M.~Pulvirenti.
\newblock On the {Derivation} of the {Hartree} {Equation} from the {N}-{Body}
  {Schrödinger} {Equation}: {Uniformity} in the {Planck} {Constant}.
\newblock {\em Journal of Functional Analysis}, 275(7):1603--1649, Oct. 2018.

\bibitem{graffi_mean-field_2003}
S.~Graffi, A.~Martinez, and M.~Pulvirenti.
\newblock Mean-{Field} {Approximation} of {Quantum} {Systems} and {Classical}
  {Limit}.
\newblock {\em Mathematical Models and Methods in Applied Sciences},
  13(01):59--73, Jan. 2003.
\newblock arXiv: math-ph/0205033.

\bibitem{gerard_homogenization_1997}
P.~Gérard, P.~A. Markowich, N.~J. Mauser, and F.~Poupaud.
\newblock Homogenization {Limits} and {Wigner} {Transforms}.
\newblock {\em Communications on Pure and Applied Mathematics}, 50(4):323--379,
  1997.

\bibitem{hauray_n-particles_2007}
M.~Hauray and P.-E. Jabin.
\newblock N-particles {Approximation} of the {Vlasov} {Equations} with
  {Singular} {Potential}.
\newblock {\em Archive for Rational Mechanics and Analysis}, 183(3):489--524,
  2007.

\bibitem{hauray_particle_2015}
M.~Hauray and P.-E. Jabin.
\newblock Particle approximation of {Vlasov} equations with singular forces:
  propagation of chaos.
\newblock {\em Annales Scientifiques de l'École Normale Supérieure.
  Quatrième Série}, 48(4):891--940, 2015.

\bibitem{hayashi_smoothing_1989}
N.~Hayashi and T.~Ozawa.
\newblock Smoothing effect for some {Schrödinger} equations.
\newblock {\em Journal of Functional Analysis}, 85(2):307--348, Aug. 1989.

\bibitem{holding_uniqueness_2018}
T.~Holding and E.~Miot.
\newblock Uniqueness and stability for the {Vlasov}-{Poisson} system with
  spatial density in {Orlicz} spaces.
\newblock In {\em Mathematical analysis in fluid mechanics—selected recent
  results}, volume 710 of {\em Contemp. {Math}.}, pages 145--162. Amer. Math.
  Soc., Providence, RI, 2018.

\bibitem{illner_global_1994}
R.~Illner, P.~F. Zweifel, and H.~Lange.
\newblock Global {Existence}, {Uniqueness} and {Asymptotic} {Behaviour} of
  {Solutions} of the {Wigner}-{Poisson} and {Schrodinger}-{Poisson} {Systems}.
\newblock {\em Mathematical Methods in the Applied Sciences}, 17(5):349--376,
  Apr. 1994.

\bibitem{jabin_mean_2016}
P.-E. Jabin and Z.~Wang.
\newblock Mean {Field} {Limit} and {Propagation} of {Chaos} for {Vlasov}
  {Systems} with {Bounded} {Forces}.
\newblock {\em Journal of Functional Analysis}, 271(12):3588--3627, Dec. 2016.

\bibitem{lazarovici_vlasov-poisson_2016}
D.~Lazarovici.
\newblock The {Vlasov}-{Poisson} {Dynamics} as the {Mean} {Field} {Limit} of
  {Extended} {Charges}.
\newblock {\em Communications in Mathematical Physics}, 347(1):271--289, 2016.

\bibitem{lazarovici_mean_2017}
D.~Lazarovici and P.~Pickl.
\newblock A {Mean} {Field} {Limit} for the {Vlasov}–{Poisson} {System}.
\newblock {\em Archive for Rational Mechanics and Analysis}, 225(3):1201--1231,
  Sept. 2017.

\bibitem{lewin_hartree_2015}
M.~Lewin and J.~Sabin.
\newblock The {Hartree} {Equation} for {Infinitely} {Many} {Particles} {I}.
  {Well}-{Posedness} {Theory}.
\newblock {\em Communications in Mathematical Physics}, 334(1):117--170, Feb.
  2015.

\bibitem{lieb_analysis_2001}
E.~H. Lieb and M.~Loss.
\newblock {\em Analysis}.
\newblock American Mathematical Society, Providence, RI, 2 edition edition,
  2001.

\bibitem{lions_sur_1993}
P.-L. Lions and T.~Paul.
\newblock Sur les mesures de {Wigner}.
\newblock {\em Revista Matemática Iberoamericana}, 9(3):553--618, 1993.

\bibitem{lions_propagation_1991}
P.~L. Lions and B.~Perthame.
\newblock Propagation of moments and regularity for the 3-dimensional
  {Vlasov}-{Poisson} system.
\newblock {\em Inventiones mathematicae}, 105(2):415--430, 1991.

\bibitem{loeper_uniqueness_2006}
G.~Loeper.
\newblock Uniqueness of the solution to the {Vlasov}–{Poisson} system with
  bounded density.
\newblock {\em Journal de Mathématiques Pures et Appliquées}, 86(1):68--79,
  July 2006.

\bibitem{markowich_classical_1993}
P.~A. Markowich and N.~J. Mauser.
\newblock The classical limit of a self-consistent quantum-vlasov equation in
  3d.
\newblock {\em Mathematical Models and Methods in Applied Sciences},
  03(01):109--124, Feb. 1993.

\bibitem{miot_uniqueness_2016}
E.~Miot.
\newblock A {Uniqueness} {Criterion} for {Unbounded} {Solutions} to the
  {Vlasov}-{Poisson} {System}.
\newblock {\em Communications in Mathematical Physics}, 346(2):469--482, 2016.

\bibitem{mitrouskas_bogoliubov_2016}
D.~Mitrouskas, S.~Petrat, and P.~Pickl.
\newblock Bogoliubov corrections and trace norm convergence for the {Hartree}
  dynamics.
\newblock {\em arXiv preprint arXiv:1609.06264}, Sept. 2016.
\newblock arXiv: 1609.06264.

\bibitem{petrat_hartree_2017}
S.~Petrat.
\newblock Hartree {Corrections} in a {Mean}-field {Limit} for {Fermions} with
  {Coulomb} {Interaction}.
\newblock {\em Journal of Physics A: Mathematical and Theoretical},
  50(24):244004, 2017.

\bibitem{petrat_new_2016}
S.~Petrat and P.~Pickl.
\newblock A {New} {Method} and a {New} {Scaling} for {Deriving} {Fermionic}
  {Mean}-{Field} {Dynamics}.
\newblock {\em Mathematical Physics, Analysis and Geometry}, 19(1):3, Mar.
  2016.

\bibitem{pickl_simple_2011}
P.~Pickl.
\newblock A {Simple} {Derivation} of {Mean} {Field} {Limits} for {Quantum}
  {Systems}.
\newblock {\em Letters in Mathematical Physics}, 97(2):151--164, Aug. 2011.

\bibitem{porta_mean_2017}
M.~Porta, S.~Rademacher, C.~Saffirio, and B.~Schlein.
\newblock Mean {Field} {Evolution} of {Fermions} with {Coulomb} {Interaction}.
\newblock {\em Journal of Statistical Physics}, 166(6):1345--1364, 2017.

\bibitem{rodnianski_quantum_2009}
I.~Rodnianski and B.~Schlein.
\newblock Quantum {Fluctuations} and {Rate} of {Convergence} {Towards} {Mean}
  {Field} {Dynamics}.
\newblock {\em Communications in Mathematical Physics}, 291(1):31--61, Oct.
  2009.
\newblock arXiv: 0711.3087.

\bibitem{saffirio_mean-field_2018}
C.~Saffirio.
\newblock Mean-field {Evolution} of {Fermions} with {Singular} {Interaction}.
\newblock {\em arXiv preprint arXiv:1801.02883}, Jan. 2018.
\newblock arXiv: 1801.02883.

\bibitem{santambrogio_optimal_2015}
F.~Santambrogio.
\newblock {\em Optimal {Transport} for {Applied} {Mathematicians}}, volume~87
  of {\em Progress in {Nonlinear} {Differential} {Equations} and {Their}
  {Applications}}.
\newblock Springer International Publishing, Cham, 2015.

\bibitem{schmeiser_vector-valued_2005}
H.-J. Schmeißer and W.~Sickel.
\newblock Vector-valued {Sobolev} {Spaces} and {Gagliardo}-{Nirenberg}
  {Inequalities}.
\newblock In {\em Nonlinear {Elliptic} and {Parabolic} {Problems}}, Progress in
  {Nonlinear} {Differential} {Equations} and {Their} {Applications}, pages
  463--472. Springer, 2005.

\bibitem{simon_trace_2005}
B.~Simon.
\newblock {\em Trace {Ideals} and {Their} {Applications}: {Second} {Edition}},
  volume 120 of {\em Mathematical {Surveys} and {Monographs}}.
\newblock American Mathematical Society, 2 edition edition, 2005.

\bibitem{thirring_quantum_1983}
W.~Thirring.
\newblock {\em Quantum {Mechanics} of {Large} {Systems}}.
\newblock Number~4 in Lehrbuch der mathematischen {Physik}. Springer, Wien,
  1983.
\newblock OCLC: 180486840.

\bibitem{villani_topics_2003}
C.~Villani.
\newblock {\em Topics in {Optimal} {Transportation}}.
\newblock American Mathematical Society, Providence, RI, Mar. 2003.

\bibitem{zhang_limit_2002}
P.~Zhang, Y.~Zheng, and N.~J. Mauser.
\newblock The {Limit} from the {Schrödinger}-{Poisson} to the
  {Vlasov}-{Poisson} {Equations} with {General} {Data} in {One} {Dimension}.
\newblock {\em Communications on Pure and Applied Mathematics}, 55(5):582--632,
  May 2002.

\end{thebibliography}

\bigskip
\bigskip
\signll

\end{document}